\documentclass[12pt]{article}

\usepackage{amssymb}
\usepackage{amsmath}
\usepackage{graphicx}
\usepackage{color}

 \newtheorem{thm}{Theorem}[section]
 \newtheorem{prop}{Proposition}[section]
 \newtheorem{lem}{Lemma}[section]

 \newtheorem{defn}[thm]{Definition}

\renewcommand{\leq}{\leqslant}
\renewcommand{\geq}{\geqslant}

\newcommand{\rbx}{\hfill{\rule{1ex}{1ex}}}

  \newcommand{\intl}{\int\limits}
  \newcommand{\suml}{\sum\limits}
  \newcommand{\liml}{\lim \limits}
  \newcommand{\tht}{\theta_\tau}
\newcommand{\gmt}{{\Gamma_\tau}}
\newcommand{\btt}{\beta_\tau}
\newcommand{\sk}{{\bf k}}

\newcommand{\diag}{\mathrm{diag}\,}

 \newcommand{\vp}{\varphi}
 \newcommand{\ve}{\varepsilon}
\newcommand{\D}{\displaystyle}

\newcommand{\ig}{\int_\Gamma}

\renewcommand{\Im}{\mathrm{Im}\,}

\newcommand{\ima}{\textrm{Im}\,}

\newcommand{\cA}{\mathcal{A}}
\newcommand{\cB}{\mathcal{B}}

\newcommand{\cG}{\mathcal{G}}

\newcommand{\cI}{\mathcal{I}}
\newcommand{\cJ}{\mathcal{J}}
\newcommand{\cK}{\mathcal{K}}
\newcommand{\cL}{\mathcal{L}}
\newcommand{\cM}{\mathcal{M}}
\newcommand{\cN}{\mathcal{N}}

\newcommand{\fL}{{\mathfrak L}}
\newcommand{\fM}{{\mathfrak M}}

\newcommand{\sC}{{\mathbb C}}

\newcommand{\sN}{{\mathbb N}}

\newcommand{\sR}{{\mathbb R}}

\newcommand{\sZ}{{\mathbb Z}}

\begin{document}

\vspace*{10mm}

\begin{center}
{\Large\textbf{Spline Galerkin methods for the double layer potential equations on contours with corners}}
\end{center}

\vspace{5mm}

\begin{center}

\textbf{Victor D. Didenko and Anh My Vu}

\vspace{2mm}


Universiti Brunei Darussalam, Bandar Seri Begawan, BE1410  Brunei

\vspace{1mm}

\texttt{diviol@gmail.com}; \texttt{anhmy7284@gmail.com}
 \end{center}

\vspace{5mm}

\textit{Dedicated to Roland Duduchava on the occasion of his seventieth
birthday}

\begin{abstract}
Spline Galerkin methods for the double layer potential equation on
contours with corners are studied. The stability of the method
depends on the invertibility of some operators $R_{\tau}$ associated
with the corner points $\tau$. The operators $R_{\tau}$ do not
depend on the shape of the contour but only on the opening angles of
the corner points $\tau$. The invertibility of these operators is
studied numerically via the stability of the method on model curves,
all corner points of which have the same opening angle. The case of
the splines of order $0,1$ and $2$ is considered. It is shown that
no opening angle located in the interval $[0.1\pi,1.9\pi]$ can cause
the instability of the method. This result is in strong contrast
with the Nystr{\"o}m method, which has four instability angles in
the interval mentioned. Numerical experiments show a good
convergence of the methods even if the right-hand side of the
equation has discontinuities located at the corner points of the
contour.
 \end{abstract}

  \vspace{6mm}


\vspace{2mm}


\section{Introduction\label{s1}}
Let $D$ be a simply connected bounded domain in $\sR^2$ with
boundary $\Gamma$, and let $n_\tau$ denote the outer normal to
$\Gamma$ at the point $\tau \in \Gamma$. It is well known that the
solution of various boundary value problems for the Laplace equation
can be reduced to solution of the integral equation
 \begin{equation}\label{eqn1}
(A_\Gamma\omega)(t) = \omega(t) + \frac{1}{\pi} \ig \omega(\tau)
\frac{d}{d n_{\tau }}\log |t-\tau |\, ds_{\tau} + (T\omega)(t)=
f(t), \quad t=x+iy \in \Gamma
\end{equation}
where $ds_{\tau}$ refer to the arc length differential and $T$ is a
compact operator. The operator
 \begin{equation*}
 V_\Gamma \omega(t):= \frac{1}{\pi} \ig \omega(\tau) \frac{d}{d n_{\tau }}\log |t-\tau
|\, ds_{\tau}
 \end{equation*}
is called the double layer potential operator and it is well known
\cite{Muskhelishvili1968} that it can be represented in the form
 $$
 V_\Gamma = \frac12(S_\Gamma + MS_\Gamma M),
 $$
where $S_\Gamma$ is the Cauchy singular integral operator on
$\Gamma$,
  $$
  (S_\Gamma x)(t): =\frac{1}{\pi i}\ig  \frac{x(\tau)\,d\tau}{\tau -t}.
  $$
and $M$ is the operator of complex conjugation,
$M\vp(t):=\overline{\vp(t)}$.

If $\Gamma$ is a smooth closed curve, then the double layer
potential operator $V_\Gamma$ is compact in the space $L^p$. This
fact essentially simplifies the stability investigation of
approximation methods for the equation \eqref{eqn1}. However, if
$\Gamma$ possesses corner points, the situation becomes more
involved. One of the simplest cases to treat is a polygonal boundary
or a boundary with polygonally shaped corners and there are a number
of works investigating approximation methods for the equation
\eqref{eqn1} on such curves \cite{Chandler:1984, Co:1985, Co:1987,
Graham:1988, Kre:1998, Kre:1990}. For a comprehensive survey, we
refer the reader to \cite{Atkinson1997, Kre:2014}.

In the present work we consider spline Galerkin methods for the
double layer potential equation \eqref{eqn1} in the case of simple
piecewise smooth curves. Such methods are often used to determine
approximate solutions of \eqref{eqn1}. However, for contours with
corners, the stability analysis of the methods is not complete. On
the other hand, it is known that for boundary integral equations the
presence of corners on the boundary may lead to extra conditions
required for the stability of the approximation method considered
\cite{DH:2011,DH:2011b,DH:2013a,VV:2014}. The aim of the present
work is twofold.  First, we obtain necessary and sufficient
conditions for the stability of spline Galerkin methods. It turns
out that stability depends on the invertibility of some operators
associated with corner points of $\Gamma$. These operators belong to
an algebra of Toeplitz operators and, at present, there is no tool
to verify their invertibility. Therefore, our second goal is to
present an approach to check the invertibility of the operators
mentioned. This approach is based on considering our approximation
methods on special model curves, and it allows us to show that
Galerkin methods for double layer potential equations on piecewise
smooth contours behave similarly to equations on smooth curves.
Thus, it was discovered that at least for the splines of degree
$0,1$ and $2$ the corresponding Galerkin method is always stable
provided that all opening angles of the corner points are located in
the interval $[0.1\pi, 1.9\pi]$. Similar results concerning the
spline Galerkin methods for the Sherman-Lauricella equation have
been recently obtained in \cite{DTV:2015}. Note that this effect is
in strong contrast with the behaviour of the Nystr\"om method which
possesses instability angles in the interval $[0.1\pi, 1.9\pi]$,
\cite{VV:2014}.

This paper is organized as follows. Section \ref{s2} is devoted to
description of spline spaces and spline Galerkin methods. Here we
also present some numerical examples illustrating the efficiency of
the method. Stability conditions are established in Section
\ref{s3}, while Section \ref{s4} deals with the numerical approach
to the search of critical angles.

\section{Spline spaces and spline Galerkin methods\label{s2}}

Let us identify each point $(x,y)$ of $\sR^2$ with the corresponding
point $z=x+iy$ in the complex plane $\sC$. By $L^2=L^2(\Gamma)$ we
denote the set of all Lebesgue measurable functions $f=f(t), t\in
\Gamma$ such that
$$
||f||=\left (\ig |f(t)|^2\,|dt|\right )^{1/2}<+\infty.
$$
By $\cM_\Gamma$ we denote the set of all corner points $\tau_0,
\tau_1, \ldots , \tau_{q-1}$ of $\Gamma$. In order to describe the
spline spaces on $\Gamma$, let us assume that this contour is
parametrized by a $1$-periodic function $\gamma:\sR\mapsto \sC$ such
that
 \begin{equation}\label{corner1}
 \tau_k=\gamma\left(\frac{k}{q}\right),\; k=0,1,\ldots ,q-1.
 \end{equation}
In addition, we also assume that the function $\gamma$ is two times
continuously differentiable on each subinterval $(k/q,(k+1)/q)$ and
\begin{equation}\label{corner2}
\left|\gamma '\left(\frac{k}{q}+0\right)\right|= \left|\gamma '
\left(\frac{k}{q}-0\right)\right|,\; k=0,1,\ldots , q-1.
 \end{equation}
For any two functions $f,g\in L^2(\sR)$, let $f*g$ denote their
convolution, i.e.
 $$
(f*g)(s):=\intl _{\sR} f(s-x)g(x)dx.
  $$
If $\chi$ is the characteristic function of the interval $[0,1)$,
then for any fixed $d\in\sN$, let
$\widehat{\phi}=\widehat{\phi}^{(d)}(s)$ be the function defined by
the recursive relation
$$\widehat{\phi}^{(d)}(s)= \begin{cases}
\chi (s) & \text{ if } d=0,\\
(\chi * \widehat{\phi}^{(d-1)})(s) &\text{ if } d>0.
\end{cases}
$$
The parametrization $\gamma$ can be now used to introduce spline
spaces on $\Gamma$. More precisely, let $n$ and $d$ be fixed
non-negative integers such that $n\geq d+1$. By $\cI(n,d)$ we denote
the set of all integers $j\in \{0,1,\ldots,n-(d+1)\}$ such that the
interval $[j/n,(j+d+1)/n]$ does not contain any point $s_k=k/q$,
$k=0,1,\ldots,q$. Let $S^d_n = S^d_n(\Gamma)$ be the set of all
linear combinations of the functions
$$
\widehat{\phi}_{nj}(t) = \widehat{\phi}^{(d)}(ns-j),\; t=\gamma(s)
\in \Gamma,\quad s\in \sR,\quad j\in \cI (n,d).
$$
For each $j\in \cI (n,d)$ set
 $$
\phi_{nj}:=\nu_d \sqrt{n} \widehat {\phi}_{nj},
 $$
where
$$
\nu_d = \left( \intl _0^{d+1} \left |\widehat{\phi}^{(d)}(s)\right
|^2 ds\right)^{-1/2}.
$$
It is easily seen that $\phi_{nj}$ are normalized
functions, i.e. $||\phi_{nj}||=1$.

\begin{figure}[tb]
  \centering
\includegraphics[height=44mm,width=61mm]{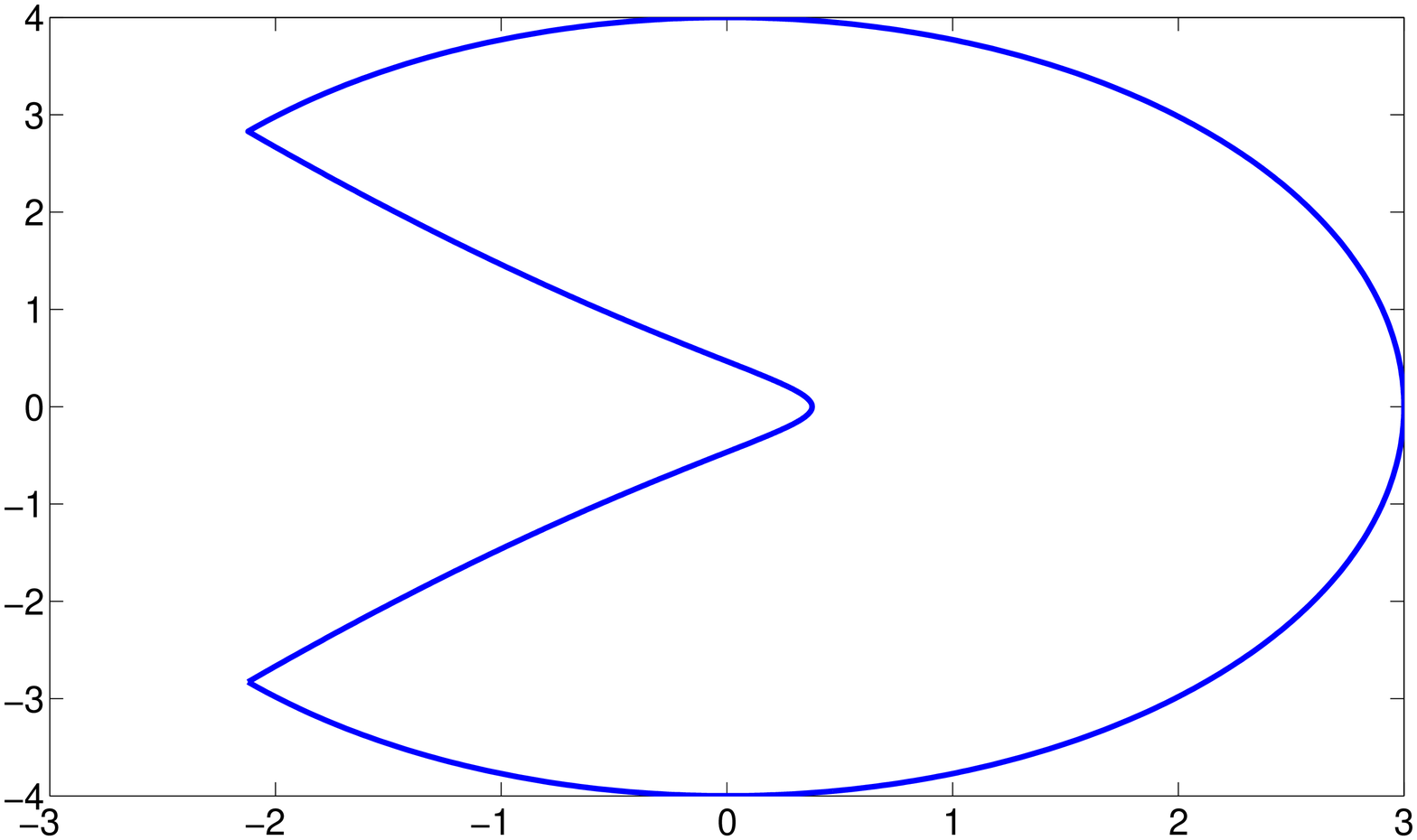}~%
\includegraphics[height=44mm,width=61mm]{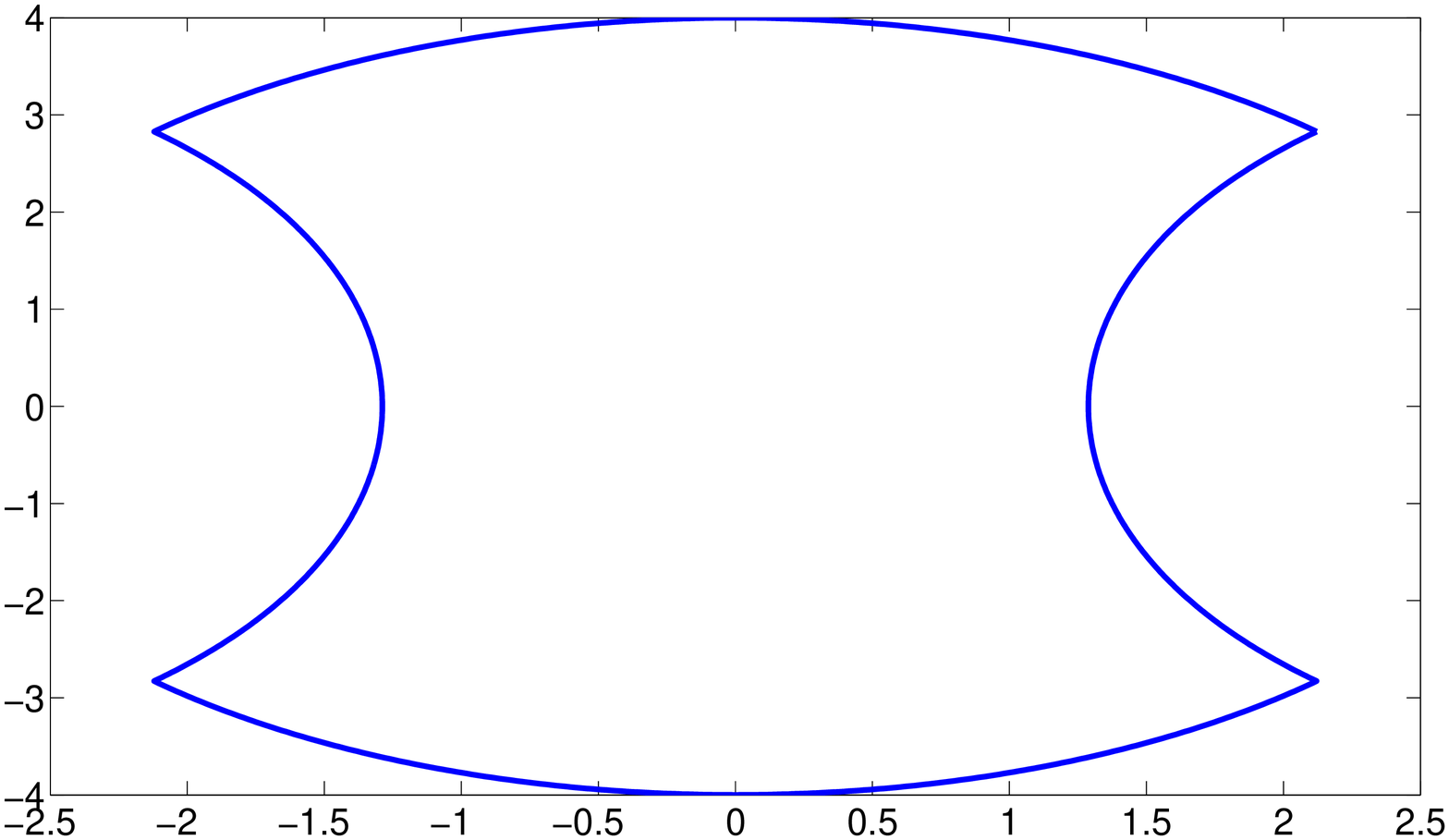}
 \caption{Left: 'pacman' curve $\fL$; Right: 'battleax' curve $\fM$.}
 \label{24cor}
\end{figure}

According to the spline Galerkin method, approximate solution
$\omega_n$ of the equation \eqref{eqn1} is sought in the form
\begin{equation}\label{glk1}
\omega _n(t) = \suml _{j \in \cI(n,d)} a_j\phi_{nj}(t)
\end{equation}
with the coefficients $a_j$ obtained from the system of linear
algebraic equations
\begin{equation}\label{glk2}
(A_\Gamma \omega_n,\phi_{nj})=(f,\phi_{nj}),\quad j \in \cI (n,d).
\end{equation}
Note that the scalar product $(\cdot, \cdot)$ is defined by
  $$
(f,g)=\intl _0^1 f(\gamma(s))\overline {g(\gamma(s))}ds.
 $$
The stability of this Galerkin method will be studied in Section
\ref{s3}. However, here we would like to illustrate the efficiency
of the method by a few examples. For simplicity, now we only
consider equations with the operator $T=0$. Although special, this
case is of the utmost importance.  It occurs when reducing boundary
value problems for partial differential equations to boundary
integral equations. In particular, we determine Galerkin solutions
of the double layer potential equation with various right-hand sides
$f$ on two curves with corners. One of these right-hand sides is
continuous on both curves, whereas two others have discontinuity
points, some of which coincide with the corners. Let us describe the
curves and right-hand sides in more details. The curves $\fL$ and
$\fM$ are obtained from the ellipse
$$
\gamma_e(s) =a\cos (2\pi s)+ib\sin (2\pi s), \quad s\in\sR,
$$
by cutting a part of it and connecting the cutting points by arcs
representing cubic Hermit interpolation polynomials in such a way
that each common point of the curve obtained becomes a corner point
satisfying the conditions \eqref{corner1}, \eqref{corner2}. In
Figure \ref{24cor}, the semi-axes of the ellipse are $a=3,b=4$. The
curve $\fL$ has two corner points obtained by cutting off the part
of the ellipse corresponding to the parameter $s \in [3/8,5/8]$. On
the other hand, two parts of the ellipse corresponding to the
parameter $s\in [3/8,5/8]\cup [7/8,9/8]$ are cut off to create the
curve $\fM$. The parametrization of the remaining parts of the
curves $\fL$ and $\fM$ is scaled and shifted so that the conditions
\eqref{corner1} and \eqref{corner2} are satisfied. Let $f_1,f_2$ and
$f_3$ be the following functions defined on the curves $\fL$ and
$\fM$,
 $$
f_1(z) = -z|z|,
 $$
 $$
 f_2(z) = \begin{cases} -1+iz &\text{ if } \ima z < 0, \\
\quad \! 1+iz &\text{ if } \ima z \geq 0, \end{cases}
 $$
and
 $$
f_3(z) = \begin{cases} -2+iz &\text{ if } \ima z < \ima z_0, \\
\quad \! 2+iz &\text{ if } \ima z \geq \ima z_0, \end{cases}
$$
where $z_0=\gamma_e (3/8)$.

In passing note that the function $f_2$ has two discontinuity points
neither of which coincides with a corner of $\fL$ or $\fM$. On the
other hand, one of the corner points of $\fL$ is a discontinuity
point for the function $f_3$, and two discontinuity points of $f_3$
are located at the corner points of $\fM$. Let $\omega_n
=\omega_n(f_j, \Gamma)$ be the Galerkin solution \eqref{glk1},
\eqref{glk2} of the double layer potential equation with right-hand
side $f_j$ considered on a curve $\Gamma$, and let
$E_n^{f_j,\Gamma}$ be the quantity
 $$
E_n^{f_j,\Gamma}=\|\omega_{2n}(f_j, \Gamma)-\omega_n(f_j,
\Gamma)\|_2/\|\omega_{2n}(f_j, \Gamma)\|_2,
 $$
which shows the rate of convergence of the approximation method
under consideration. The Table \ref{convergence} illustrates how the
spline Galerkin method with $d=0$ performs for the curves $\fL$ and
$\fM$ and for the right-hand sides $f_1,f_2$ and $f_3$.

\begin{figure}[tb]
   \centering
\includegraphics[height=36mm,width=60mm]{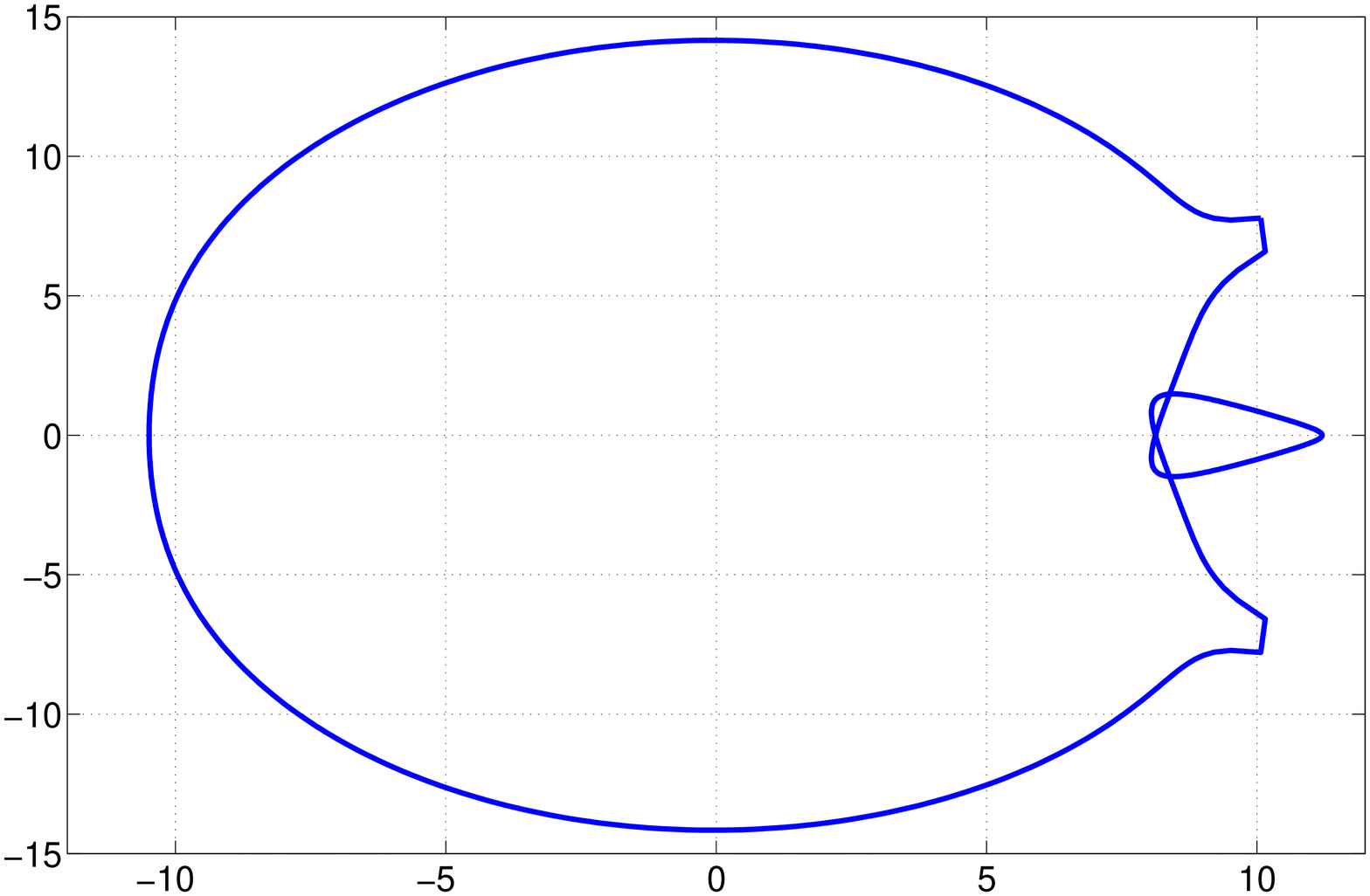}~%
 \includegraphics[height=36mm,width=60mm]{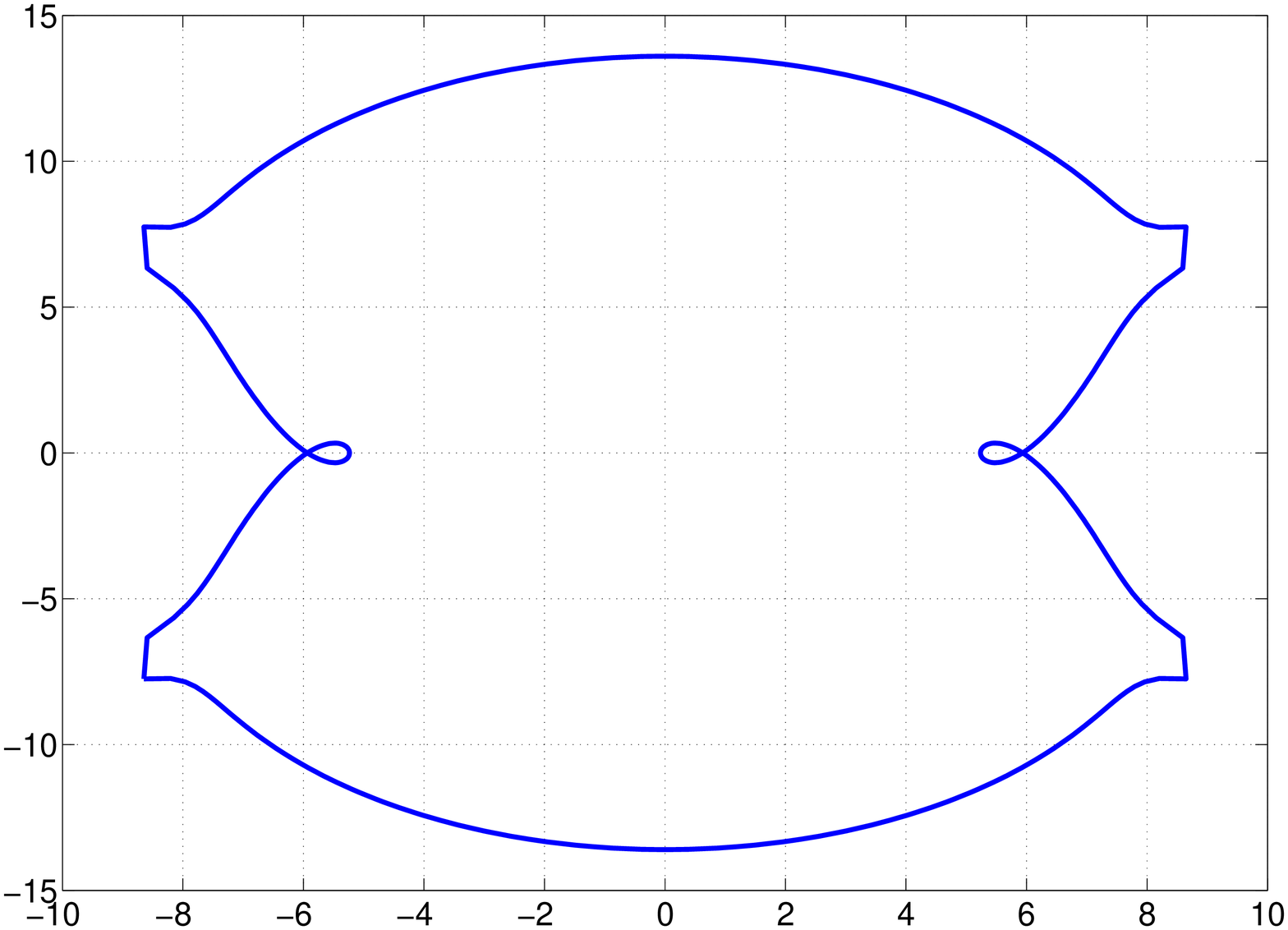}
 \includegraphics[height=36mm,width=60mm]{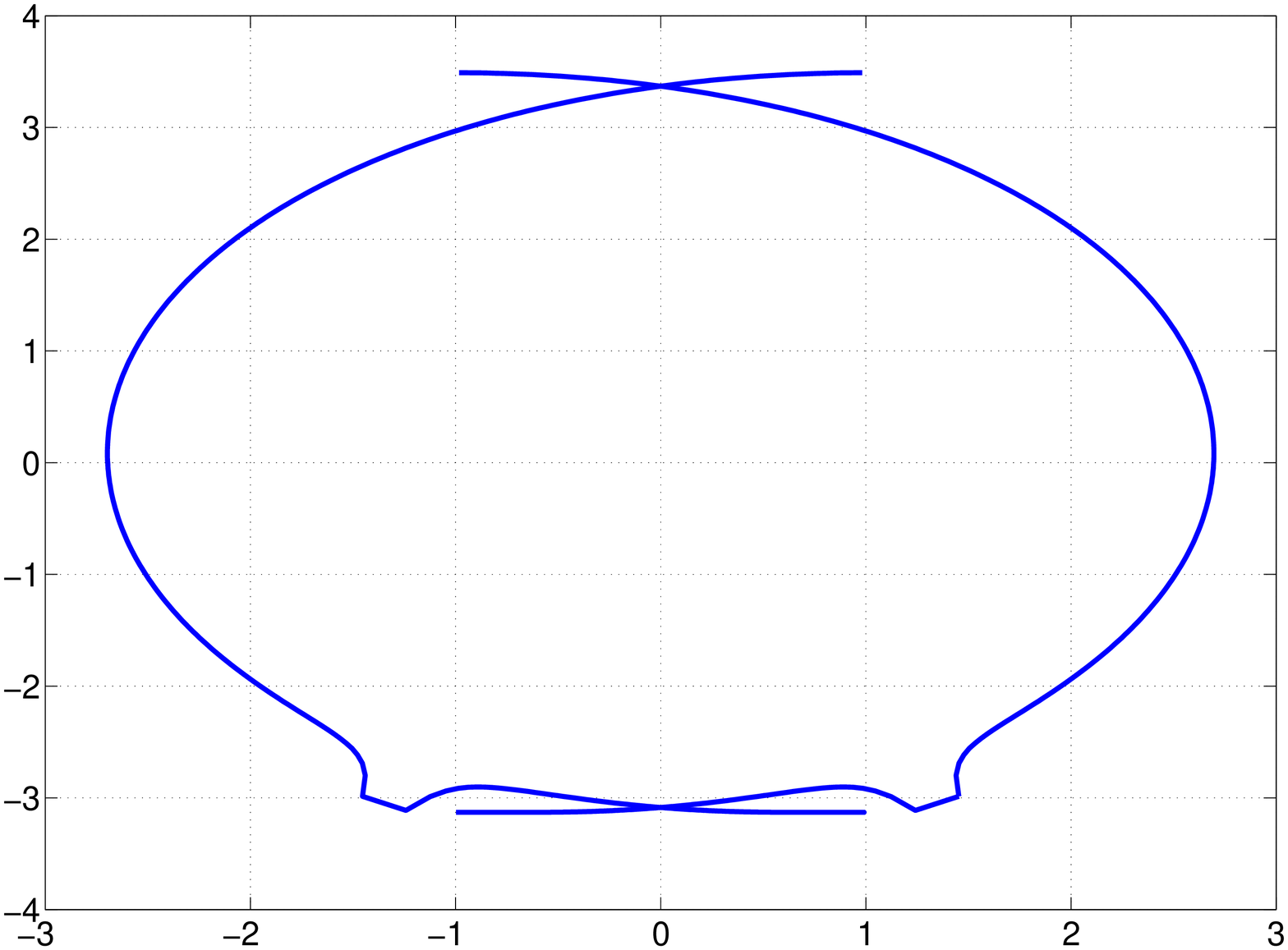}~%
 \includegraphics[height=36mm,width=60mm]{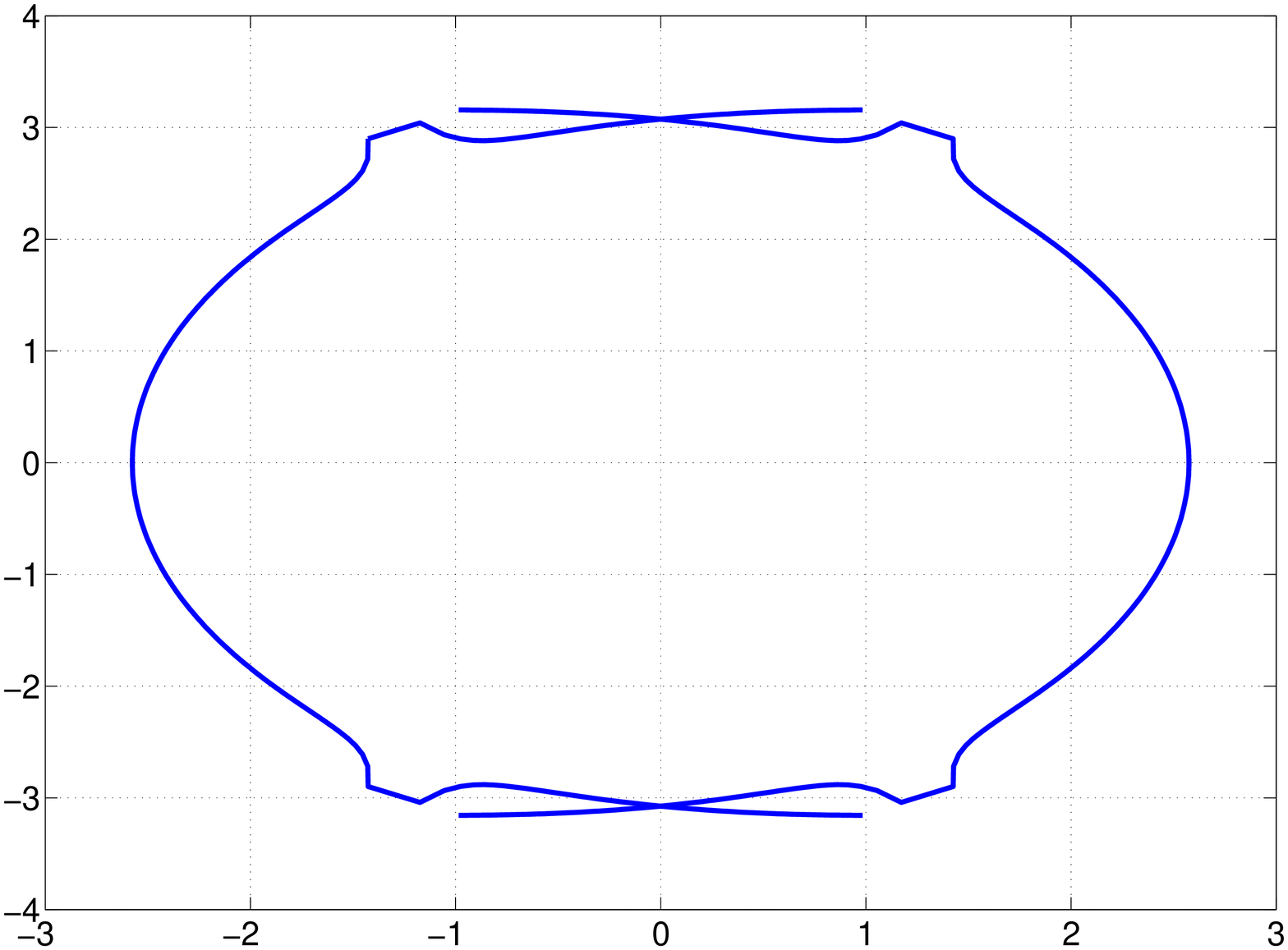}
 \includegraphics[height=36mm,width=60mm]{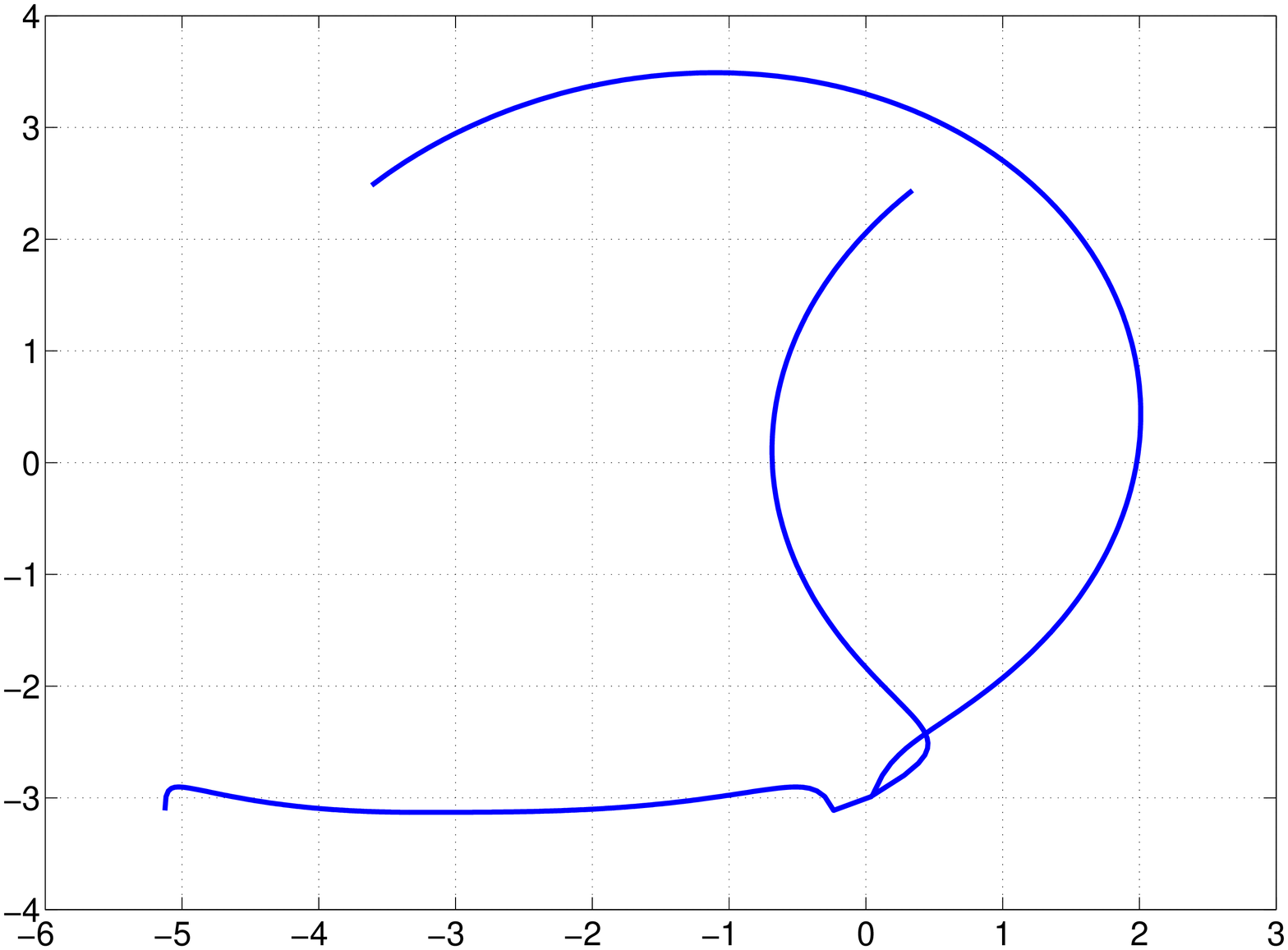}~%
 \includegraphics[height=36mm,width=60mm]{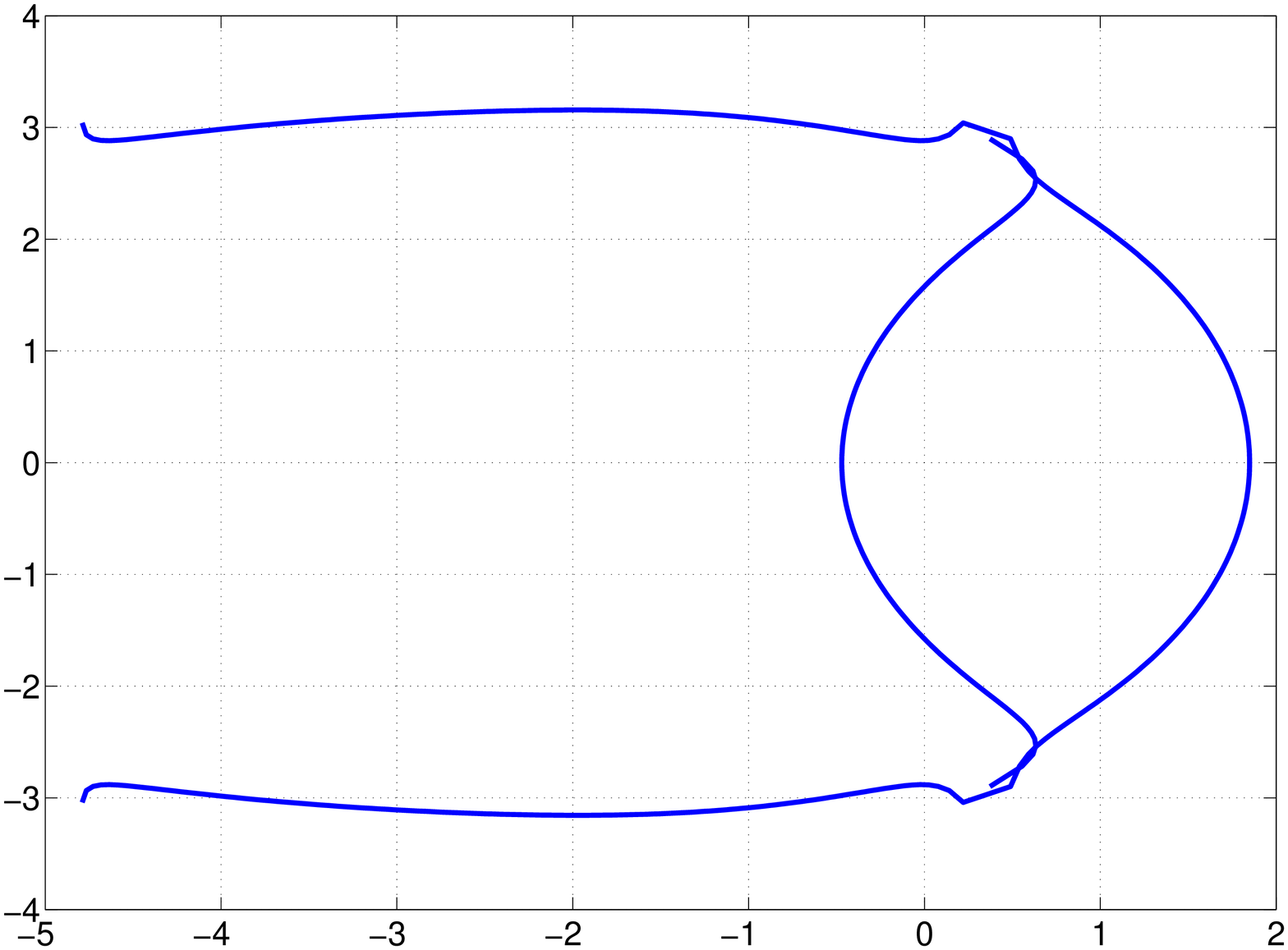}
 \caption{Spline Galerkin solution of the double layer potential equation
in the case $n=1024$. Left: 'pacman' curve $\fL$; Right: 'battleax'
curve $\fM$. First row: with r.-h. s. $f_1(z)$; Second row: with
r.-h. s. $f_2(z)$; Third row: with r.-h. s.$f_3(z)$.}
 \label{app_sol}
\end{figure}

Note that the integrals in the scalar products
$(A_\Gamma\omega_n,\phi_{nj}),\; j\in \cI(n,d)$ are approximated by
the Gauss-Legendre quadrature formula with quadrature points
coinciding with the zeros of the Legendre polynomial of degree $24$
on the canonical interval $[-1,1]$ scaled and shifted to the
intervals $[j/n,(j+d+1)/n]$. More precisely, we employ the formula
\begin{equation}
(A_\Gamma\omega_n,\phi_{nj}) =\intl _0^1 A_\Gamma\omega_n
(\gamma(s)) \overline{\phi_{nj}(\gamma(s))}ds \approx \suml
_{k=1}^{24}w_k A_\Gamma\omega_n(\gamma(s_k)) \phi_{nj}(\gamma(s_k)),
\end{equation}
where $w_k,s_k$ are weights and Gauss-Legendre points on the
interval $[j/n,(j+d+1)/n]$. Composite Gauss-Legendre quadrature
 is also used in approximation of the
integral operators of $A_\Gamma\omega_n(\gamma(s_k))$, cf.
\cite{DH:2011}. Thus we employ the quadrature formula
$$
\begin{aligned}
\ig k(t,\tau)x(\tau) d\tau &=\intl _0^1 k(\gamma(\sigma),\gamma(s))x(\gamma(s))\gamma '(s) ds \\
                & \approx \suml _{l=0}^{m-1} \suml_{p=0}^{r-1} w_pk(\gamma(\sigma),
                \gamma(s_{lp}))x(\gamma(s_{lp}))\gamma '(s_{lp})/m
\end{aligned}
$$
where $m=40$, $r=24$, $s_{lp}=(l+\ve_p)/m$ and $w_p$ and $\ve_p$
are, respectively, the Gauss-Legendre weights and Gauss-Legendre
nodes scaled and shifted to the interval $[0,1]$. For the discrete
norm used in error evaluation, we set $h_0=1/128$, $h=10^{-3}$ and
choose the meshes
\begin{equation*} U_1:=(0+h_0:h:0.5-h_0) \bigcup
(0.5+h_0:h:1-h_0)
 \end{equation*}
 and
 \begin{align*}
U_2:=&(0\!+\!h_0:h:0.25\!-\!h_0)\bigcup
(0.25\!+\!h_0:h:0.5\!-\!h_0)\\
 &\bigcup (0.5+h_0:h:0.75-h_0) \bigcup
(0.75+h_0:h:1-h_0)
 \end{align*}
due to the fact that the curves $\fL$ and $\fM$ have two and four
corner points, respectively, cf. Condition \ref{corner1}. In the
graphs of Figure \ref{app_sol}, jumps appear when the corner points
of $\fM$ and the discontinuity points of the right-hand side $f_3$
coincide. At the same time, it is quite remarkable that the
condition numbers of the methods are relatively small. For the
interval considered, they do not exceed $10$ and $5$ for the curve
$\fL$ and $\fM$, respectively.

\begin{table}
\begin{center}
\begin{tabular}{|c||c|c||c|c||c|c|}
\hline
  n &    $E_n^{(f_1,\fL)}$ & $E_n^{(f_1,\fM)}$ & $E_n^{(f_2,\fL)}$& $E_n^{(f_2,\fM)}$ & $E_n^{(f_3,\fL)}$&$E_n^{(f_3,\fM)}$ \\
\hline
$128$ & $0.0257$ &          $0.0279 $ &               $0.0248$      &$0.0261$        &$0.0445$              &$ 0.0383$ \\
\hline
$256$ & $0.0129$&$           0.0147$ &                $0.0125$     &$0.0141$         &$0.0286$             &$0.0230$\\
\hline
$512$ &$0.0054 $&$          0.0073$  &                 $0.0052$    &$0.0070$         &$0.0186$            &$0.0153 $\\
\hline
\end{tabular}
\caption{Convergence of the spline Galerkin method, $d=0$.}
\label{convergence}
\end{center}
\end{table}
Let us also point out that the results presented in Table
\ref{convergence} are comparable with the convergence rates of the
spline Galerkin methods for the Sherman-Lauricella \cite{Xu:2014}
and Muskhelishvili \cite{VVe:2007} equations on smooth curves. These
estimates can still be improved if one uses a more accurate
approximation of the integrals arising in the Galerkin method
\cite{Hels:2011, HO:2008}. Nevertheless, the approximate solutions
presented in Figure \ref{app_sol} demonstrate a good accuracy. We
also computed Galerkin method solutions of the double layer
potential equation with the right-hand sides and curves from
\cite{VV:2014}. Although these results are not reported here, there
is a good correlation with approximate solutions of \cite{VV:2014}
obtained by the Nystr\"om method.

\section{Local operators and stability of the spline Galerkin method\label{s3}}
Let us briefly describe the approach we use in the study of the
stability of the Galerkin method. For more details, we refer the
reader to \cite{DS:2008, Prossdorf1991, RSS:2011}. Let $P_n$ denote
the orthogonal projection from $L^2(\Gamma)$ onto $S_n^d(\Gamma)$.
The spline Galerkin method \eqref{glk2} can be rewritten as
\begin{equation}\label{glk3}
P_n A_\Gamma P_n \omega_n=P_nf,\; n\in \sN.
\end{equation}
 \begin{defn}
The approximation sequence $(P_nA_\Gamma P_n)$ is said to be stable
if there exists $n_0\in \sN$ and a constant $C>0$ such that for all
$n\geq n_0$ the operators $P_n A_\Gamma P_n:S^d_n(\Gamma) \mapsto
S^d_n(\Gamma)$ are invertible and $\|(P_n A_\Gamma P_n)^{-1}P_n\|
\leq C$.
 \end{defn}
Let $\cA^\Gamma$ denote the set of all bounded sequences of bounded
linear operator $A_n:\Im P_n \mapsto \Im P_n$ such that there exist
strong limits
$$
s-\lim A_nP_n =A,\quad s-\lim (A_nP_n)^*P_n=A^*.
$$
Moreover, let $\cK (L^2(\Gamma))$ denote the ideal of all compact
operators in $\cL(L^2(\Gamma))$, and let $\cG \subset \cA^\Gamma$ be
the set of sequences which converge uniformly to zero. Recall that
 the sequence of orthogonal projection $(P_n)$ in $L^2(\Gamma)$
converges strongly to identity operator and $P_n^*=P_n$. It follows
that
$$
s-\liml_{n\to \infty} P_nA_\Gamma P_n=A_\Gamma, \quad s-\liml_{n\to
\infty} (P_nA_\Gamma P_n)^*P_n=A^*_\Gamma.
$$
It is well known\cite{DS:2008, Prossdorf1991, RSS:2011} that the set
of sequences
$$
\cJ^\Gamma =\{(A_n)\in \cA^\Gamma: A_n=P_nKP_n + G_n,\; K\in
\cK(L^2(\Gamma)),\; (G_n)\in \cG\}
$$
forms a close two-sided ideal of $\cA^\Gamma$.
\begin{prop}[{cf. \cite[Proposition 1.6.4]{DS:2008}}]\label{prop1}
The sequence $(P_nA_\Gamma P_n)$ is
stable if and only if the operator $A_\Gamma\in \cL(L^2(\Gamma))$
and the coset $(P_nA_\Gamma P_n)+\cJ^\Gamma \in \cL(\cA^\Gamma /
\cJ^\Gamma)$ are invertible.
\end{prop}
Recall that both Fredholm properties and invertibility of the
operator $A_\Gamma$ in various spaces have been studied in
literature \cite{Atkinson1997,Co:1988,Muskhelishvili1968,Ver:1984}.
Therefore, our main task here is to investigate the behaviour of the
coset $(P_nA_\Gamma P_n)+\cJ^\Gamma$. Note that it is more
convenient to consider this coset as an element of a smaller
algebra.

Thus let $\cB ^\Gamma$ denote the smallest closed $C^*$-subalgebra
of $\cA^\Gamma$ which contains the sequences $(P_nM S_\Gamma M
P_n),(P_n S_\Gamma P_n)$, all sequences $(G_n)\in \cG$, and all
sequences $(P_nfP_n)$ with $f\in C_{\sR}(\Gamma)$. It follows from
\cite{Prossdorf1991,RSS:2011} that $\cJ^\Gamma \subset \cB^\Gamma$
and $(P_nA_\Gamma P_n) \in \cB^\Gamma$. Therefore, $\cB^\Gamma /
\cJ^\Gamma$ is a $C^*$-subalgebra of $\cA^\Gamma / \cJ^\Gamma$,
hence the coset $(P_nA_\Gamma P_n)+\cJ^\Gamma$ is invertible in
$\cA^\Gamma / \cJ^\Gamma$ if and only if it is invertible in
$\cB^\Gamma / \cJ^\Gamma$. However, the invertibility of the coset
$(P_nA_\Gamma P_n)+\cJ^\Gamma$ in the quotient algebra $\cB^\Gamma /
\cJ^\Gamma$ can be showed by a local principle. Thus, with each
point $\tau \in \Gamma$, we associate a curve $\Gamma _\tau$ as
follows. Let $\theta _\tau\in (0,2\pi)$ be the angle between the
right and left semi-tangents to $\Gamma$ at the point $\tau$.
Further,  let $\beta _\tau\in (0, 2\pi)$ be the angle between the
real axis $\sR$ and the right semi-tangent to $\Gamma$ at the same
point $\tau$. Let $\Gamma_\tau$ be the curve defined by
$$
\gmt := e^{i(\btt + \tht)}\sR^+_-\, \bigcup \, e^{i\btt}\sR_+^+
$$
where $\sR^+_-$ and $\sR_+^+$ are positive semi-axes directed to and
away from zero, respectively. On the curve $\Gamma_\tau$ consider
the corresponding double layer potential operator $A_\gmt =
I+V_\gmt$. Moreover, let
$$
\widetilde{\phi}_{nj}(t) = \begin{cases}
\begin{cases}
\widehat{\phi}^{(d)}(ns-j) & \text{ if } t=e^{i\btt}s \\
0 & \text{ otherwise }
\end{cases}  \quad j \geq 0, \\
\begin{cases}
\widehat{\phi}^{(d)}(ns-j+d) & \text{ if } t=e^{i(\btt + \tht)}s \\
0 & \text{ otherwise }
\end{cases}  \quad j <0 \\
\end{cases}.
$$
By $S_n^d(\gmt)$  we denote the smallest closed subspace of
$L^2(\gmt)$ which contains all functions $\widetilde{\phi}_{nj}$,
$j\in \sZ$. Correspondingly, $S_n^d(\sR^+)$ is the smallest subspace
of $L^2(\sR^+)$ containing all functions $\widetilde{\phi}_{nj}$,
$j\geq 0$ for $\btt =0$. In addition, let $P_n^\tau$ and  $P_n^+$ be
the orthogonal projections of $L^2(\gmt)$ onto $S_n^d(\gmt)$ and
$L^2(\sR^+)$ onto $S_n^d(\sR^+)$, respectively. Now algebra $\cB
^\gmt$ and its ideal $\cJ ^\gmt$ can be defined analogously to the
construction of $\cB^\Gamma$ and $\cJ^\Gamma$.

Similarly to Proposition \ref{prop1}, one can formulate the
following result.
\begin{lem} \label{coro1}
The sequence $(P_n^\tau A_\gmt P_n^\tau) \in \cB^\gmt$ is stable if
and only if the operator $A_\gmt$ is invertible and the coset
$(P_n^\tau A_\gmt P_n^\tau) + \cJ^\gmt$ is invertible in the
quotient algebra $\cB^\gmt / \cJ^\gmt$.
\end{lem}

Note that the invertibility of the operator $A_{\Gamma_\tau}$ can be
studied quite easily. It turns out that this operator is
isometrically isomorphic to a block Mellin operator (see relations
\eqref{iso}-\eqref{smb} below). The invertibility of block Mellin
operators depends on the invertibility of their symbols in an
appropriate function algebra and is well understood \cite{Du:1979,
Simonenko1986}. It follows that the operator $A_{\Gamma_\tau}$ is
always invertible. Therefore, the coset $(P_n^\tau A_\gmt
P_n^\tau)+\cJ^\gmt$ is invertible in the corresponding quotient
algebra if and only if the sequence $(P_n^\tau A_\gmt P_n^\tau)$ is
stable. Let us now consider the stability problem in more detail. By
$L^2_2(\sR^+)$ we denote the product of two copies of $L^2(\sR^+)$
provided with the norm
$$ \|(\varphi_1,\varphi_2)^T\|_{L^2_2(\sR^+)} = \left(\|\varphi_1\|_{L^2(\sR^+)}^2+\|\varphi_2\|_{L^2(\sR^+)}^2\right)^{1/2},$$
and let $\eta:L^2(\gmt) \mapsto L^2_2(\sR^+)$ be the mapping defined
by
$$ \eta(f)=\left(f(se^{i(\btt+\tht)}),f(se^{i\btt})\right)^T.$$
This isometry generates an isometric algebra isomorphism $\Psi:
\cL(L^2(\gmt)) \mapsto \cL(L^2_2(\sR^+))$ defined by
\begin{equation}\label{iso}
\Psi(A)=\eta A \eta ^{-1}.
\end{equation}
In particular, for the operators $P_n^\tau$, $I$ and $V_\gmt$, one
has
$$
\Psi(P_n^\tau)=\diag(P_n^+,P_n^+)
$$
and
 \begin{equation}\label{psiA}
  \Psi(I) = \left [ \begin{matrix} I & 0  \\
  0 & I \end{matrix}\right ], \quad
  \Psi(V_\gmt) = \left [ \begin{matrix} 0 & \cN _{\tht}  \\
  \cN _{\tht} & 0 \end{matrix}\right ],
 \end{equation}
where $\cN _{\theta}$ is the Mellin convolution operator defined by
   \begin{equation}\label{dlp}
(\cN _\theta (\varphi))(\sigma) = \frac{1}{2\pi i}\int _0^{+\infty}
  \left(\frac{1}{s-\sigma e^{i\theta}}-\frac{1}{s-\sigma e^{-i\theta}}\right)\varphi (s)\, ds.
 \end{equation}
The operator $\cN _{\theta}$ can also be written in another form
reflecting its Mellin structure-viz.,
\begin{equation}\label{ntheta}
\cN _{\theta}(\varphi)(\sigma) = \intl_0^{+\infty} \sk_\theta \left
(\frac{\sigma}{s}\right)\varphi (s) \, \frac{ds}{s}
\end{equation}
where
  \begin{equation}\label{smb}
\sk_\theta=\sk_\theta (u) = \frac{1}{2\pi}\, \frac{u\sin
\theta}{|1-ue^{i\theta}|^2}.
 \end{equation}
Thus
 \begin{equation}\label{mell}\Psi(A_{\Gamma_\tau}) = \left [ \begin{matrix} I & \cN _{\tht}  \\
  \cN _{\tht} & I \end{matrix}\right ],
 \end{equation}
and an immediate consequence of the isomorphism \eqref{iso} is that
the sequence $(P_n^\tau A_\gmt P_n^\tau)$ is stable if and only if
so is the sequence $(\Psi(P_n^\tau A_\gmt P_n^\tau))$. On the other
hand, the study of the stability of the sequences $(\Psi(P_n^\tau
A_\gmt P_n^\tau))$, $\tau \in\Gamma$ can be reduced to the study of
two main cases related to the nature of the points $\tau\in \Gamma$.
Thus if $\tau \not \in \cM_\Gamma$, then $\tht =\pi$ so that the
operator $\cN_\pi=0$ and $\Psi(P_n^\tau A_\gmt P_n^\tau)$ is just
the diagonal sequence $\diag(P_n^+,P_n^+)$ which is obviously
stable. Therefore the corresponding coset $(P_n^\tau) + \cJ^\gmt \in
\cB^{\Gamma_\tau}/\cJ^{\Gamma_\tau}$ is invertible. Consider now the
case where $\tau$ is a corner point of $\Gamma$, and $\tht\in
(0,2\pi)$ is the opening angle of this corner. By $l^2$ we denote
the set of sequences of complex numbers $(\xi _k)_{k=0}^{+\infty}$
such that
$$
\suml _{k=0}^{\infty} |\xi _k|^2 < \infty.
$$
Moreover, let $\Lambda _n$ be the operator acting from
$S^d_n(\sR^+)$ into $l^2$ and defined by
 $$
 \Lambda _n\left(\suml _{j=0}^{\infty} \xi _j \widetilde{\phi}_{nj}\right) = (\xi _0, \xi _1, \ldots).
 $$
The operators $\Lambda _n$ are continuously invertible and there is
a constant $m$ such that
 $$
 ||\Lambda_n||\,||\Lambda_{-n}|| \leq m \quad \text{for all}\quad
 n=1,2, \ldots,
 $$
where $\Lambda_{-n}:=\Lambda_n^{-1}$, \cite{DeBoor1978}.  This
implies that the sequence $(\Psi (P_n^\tau A_\gmt P_n^\tau))$ is
stable if and only if the sequence $(R^\tau_n)$,
$$
R^\tau_n := \diag (\Lambda_n,\Lambda_n)\, \Psi(P_n^\tau A_\gmt
P_n^\tau)\, \diag(\Lambda_{-n},\Lambda_{-n}): l^2\times l^2 \mapsto
l^2\times l^2
$$
is stable.
\begin{lem}\label{lem2}
The sequence $(P_n^\tau A_\gmt P_n^\tau)$  is stable if and only if
the operator $R_\tau:=R^\tau_1$ is invertible.
\end{lem}
\textbf{Proof.}  According to the above considerations, the sequence $(P_n^\tau A_\gmt
P_n^\tau)$ is stable if and only if so is the sequence $(R_n^\tau)$.
Note that the operators $R_n^\tau$ have the form
$$
R_n^\tau=\left[\begin{matrix} I & A_{12} \\ A_{21} & I
\end{matrix}\right]
$$
where $A_{21} = A_{12}=\Lambda_n P_n^+ \cN_{\theta_\tau}P_n^+
\Lambda_{-n}$.

 Consider now the operators $\Lambda_n P_n^+ \cN_{\theta_\tau}P_n^+
\Lambda_{-n}$, $n\in \sN$. Let $J_n:S_n^d(\sR^+)\to S_1^d(\sR^+)$ be
the operator defined by
 $$
(J_n \phi)(s):=\phi(s/n).
 $$
Then $J_n\widetilde{ \phi}_{jn}=\widetilde{\phi}_{j1}$ and
$\Lambda_n=\Lambda_1 J_n$. This implies the relation
 $$
\Lambda_n P_n^+ \cN_{\theta_\tau}P_n^+ \Lambda_{-n} = \Lambda_1 J_n
P_n^+ \cN_{\theta_\tau}P_n^+ J_n^{-1} \Lambda_{-1}= \Lambda_1 P_1^+
\cN_{\theta_\tau} P_1^+ \Lambda_{-1},
 $$
so that $(R_n^\tau)$ is a constant sequence. Therefore it is stable
if and only if the operator $R_\tau=R_1^\tau$, is invertible.
 \rbx

Using the above results, one can obtain a stability criterion for
the spline Galerkin method.
\begin{thm}\label{the1}
If operator $A_\Gamma$ is invertible, then the spline Galerkin
method \eqref{glk3} is stable if and only if all the operators
$R_\tau: l^2\times l^2 \mapsto l^2\times l^2$,  $\tau \in
\cM_\Gamma$ are invertible.
\end{thm}
\textbf{Proof.}
It follows from Proposition \ref{prop1} that the spline Galerkin
method is stable if and only if the coset $(P_nA_\Gamma
P_n)+\cJ^\Gamma\in\cB^\Gamma / \cJ^\Gamma$ is invertible. By Allan's
local principle \cite{Allan1968, DS:2008, RSS:2011}, this coset is
invertible if and only if so are all the cosets $(P_n^\tau A_\gmt
P_n^\tau)+\cJ^\gmt$, $\tau \in \Gamma$. However, as we already know,
for $\tau\notin \cM_\Gamma$ the coset $(P_n^\tau A_\gmt
P_n^\tau)+\cJ^\gmt$ is always invertible in the corresponding
quotient-algebra $\cB^{\Gamma_\tau}/\cJ^{\Gamma_\tau}$. On the other
hand, if $\tau\in \cM_\Gamma$, then by the Lemma  \ref{lem2} the
coset $(P_n^\tau A_\gmt P_n^\tau)+\cJ^\gmt$ is invertible if and
only if so is the corresponding operator $R_\tau$, and the proof is
completed.
\rbx

\section{Numerical approach to the invertibility of local operators\label{s4}}

Theorem \ref{the1} shows that the stability of the spline Galerkin
method depends on the invertibility of the operators $R_\tau$,
 $\tau\in \cM_\Gamma$. However, these operators belong to an
algebra of Toeplitz operators generated by piecewise continuous
matrix functions and at present there is no analytic tool to check
their invertibility. On the other hand, a numerical approach to such
a kind of problem has been proposed in \cite{DH:2011b, DH:2013a}.
Thus one can consider stability of an approximation method on curves
having corner points with the same opening angle. If this is the
case, the stability of the corresponding method depends on the
operator itself and on only one additional operator $R_\tau$. More
precisely, the following result is true.
  \begin{prop}\label{cor2}
If $\Gamma$ is a piecewise smooth curve such that all corners
$\tau\in \cM_\Gamma$ have the same opening angle, then
 \begin{enumerate}
   \item For any $\tau_1, \tau_2 \in \cM_\Gamma$ one has
     $R_{\tau_1}=R_{\tau_2}$.
\item The operator $R_{\tau_1}$ is invertible if and only if the
spline Galerkin method $(P_n (I+V_{\Gamma})P_n)$ is stable.
            \end{enumerate}
  \end{prop}
 \textbf{Proof.} 
This result is an immediate consequence of Theorem \ref{the1}. One
only has to take into account that if $\Gamma$ satisfies the
conditions stated, then the corresponding operator $I+V_\Gamma$ is
invertible on the space $L^2(\Gamma)$, \cite{Ver:1984}.
 \rbx

\begin{figure}[!tb]
   \centering
 \includegraphics[width=85mm]{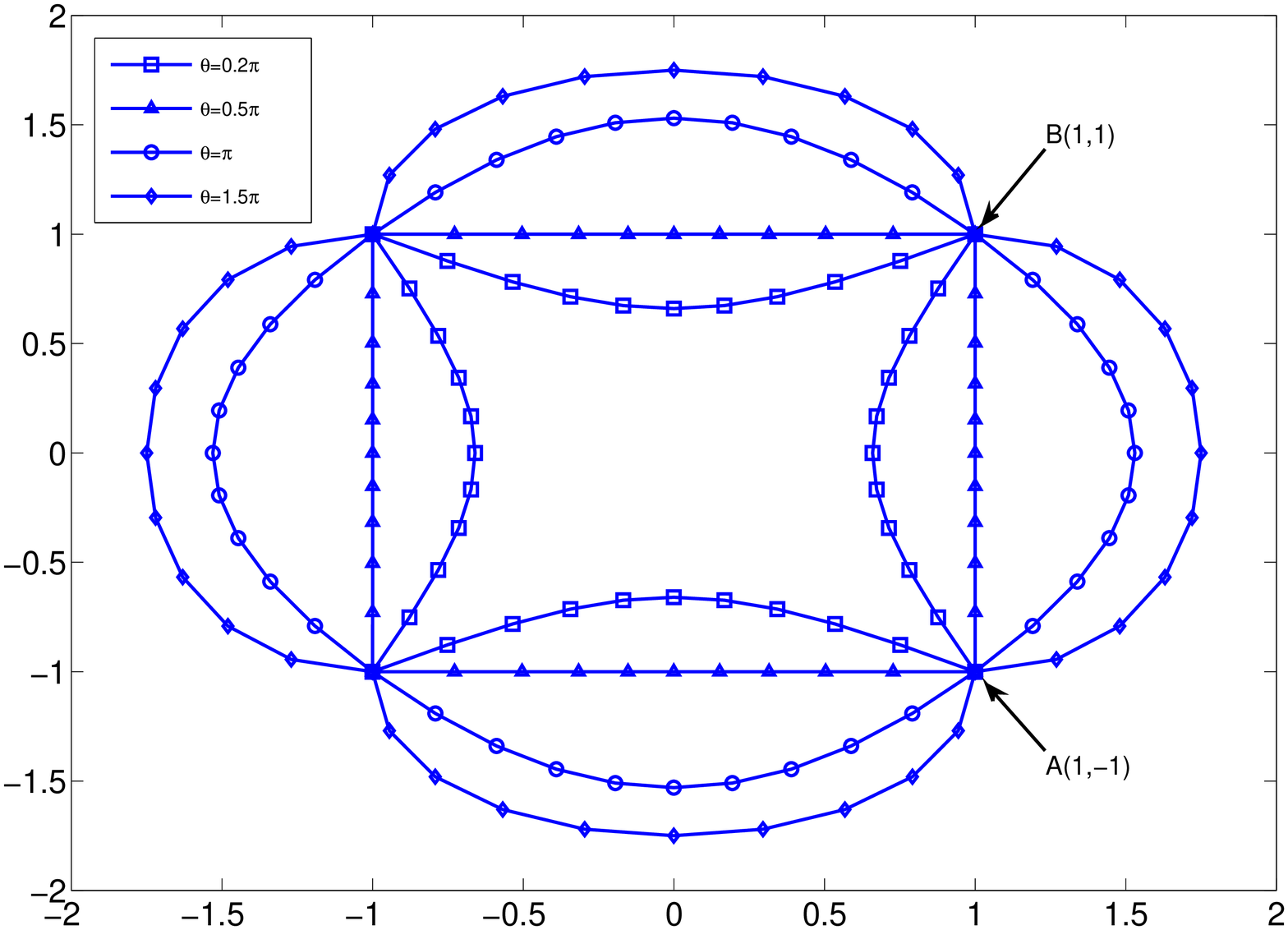}
  \caption{Curves $\fL_4(\theta)$ for various $\theta$.}
 \label{four_ang}
\end{figure}

Thus in order to detect critical angles, i.e. the opening angles
$\tht$ where the operators $R_\tau$ are not invertible, one can
compute the condition numbers of the method on families
 of special contours $\fL(\theta)$, $\theta\in (0,2\pi)$ all corner
points of which have the same opening angle. As a result, at any
critical point of the method, the graph representing the condition
numbers has to have an "infinite" peak regardless of the family of
the curves used. In this paper we employ the curves $\fL_1(\theta),
\fL_2(\theta)$, proposed in \cite{DH:2011b, DH:2013a}, which have
one and two corner points, respectively, together with a new
$4$-corner curve $\fL_4(\theta)$. The curves $\fL_1(\theta),
\fL_2(\theta)$ have the following parametrizations
\begin{align*}
&\fL_1(\theta): \gamma_1(s)=\sin(\pi s)e^{i\theta(s-0.5)},\; 0\leq s \leq 1;\\
&\fL_2(\theta): \gamma_2(s)=\left \{
\begin{array}{ll}
\D -\frac{1}{2}\cot(\theta/2)+\frac{1}{2\sin(\theta/2)}\, e^{\imath
\theta(2s-0.5)} &\quad 0\leq s<0.5,\\[1.5ex]
\D \phantom{-}\frac{1}{2}\cot(\theta/2)-\frac{1}{2\sin(\theta/2)}\,
e^{\imath \theta(2s-1.5)} &\quad 0.5\leq s< 1.
\end{array}
 \right .
\end{align*}
The $4$-corner curve $\fL_4(\theta)$ is constructed as follows.
First, connect the two points $A=(1-i)$ and $B=(1+i)$  by an arc
representing a Hermit interpolation polynomial such that
  $$
\widehat{\left(\overrightarrow{OA},
\overrightarrow{t_A}\right)}=-\theta /2,\;
\widehat{\left(\overrightarrow{OB},
\overrightarrow{t_B}\right)}=\theta/2,
 $$
where $O$ denotes the origin, $\overrightarrow{t_A}$ and
$\overrightarrow{t_B}$ are, respectively, the tangential vectors at
the points $A$ and $B$, and $\widehat{(\overrightarrow{t_1},
\overrightarrow{t_2})}$ is the angle measured from
$\overrightarrow{t_1}$ to $\overrightarrow{t_2}$ in the
counterclockwise direction (see Figure~\ref{four_ang}). Further,
rotate the arc obtained around the origin by angles $0.5\pi,\,\pi$
and $1.5 \pi$. Some curves from this family are presented in
Figure~\ref{four_ang}. Note that the Hermit interpolation polynomial
used has the parametrization
$$
P_3(s) = 1-i+as+(2i-a)s^2+(a+b-4i)s^2(s-1),\quad 0 \leq s\leq 1,
$$
where
 $$
a=3\sin\left(\frac{3\pi}{4}+\frac{\theta}{2}\right)+3i\cos\left(\frac{3\pi}{4}+\frac{\theta}{2}\right),\quad
b=3\sin\left(\frac{\pi}{4}-\frac{\theta}{2}\right)+3i\cos\left(\frac{\pi}{4}-\frac{\theta}{2}\right).
 $$

\begin{figure}[!tb]
   \centering
\includegraphics[width=43mm,height=30mm]{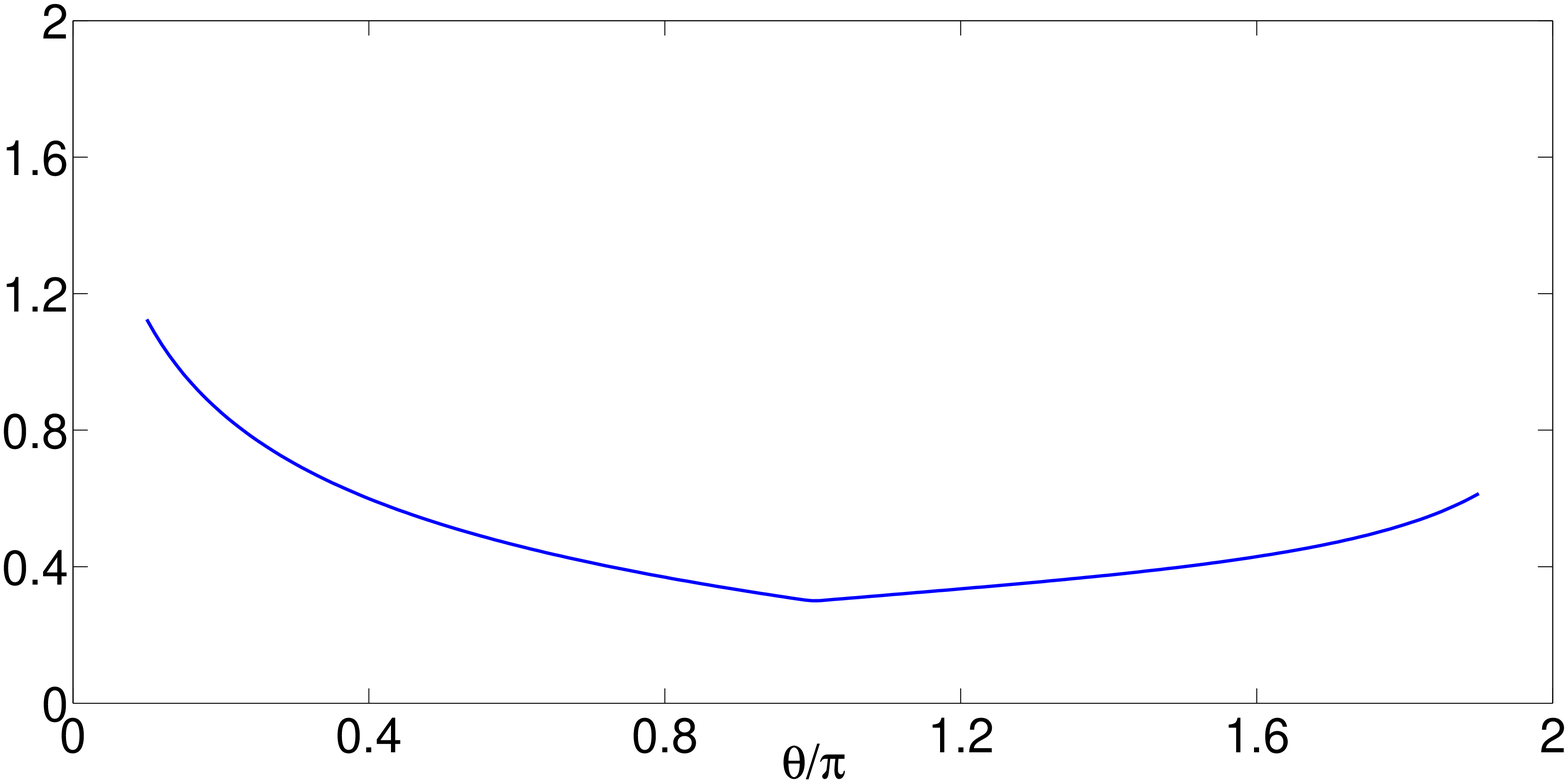}\hspace{-4mm}
\includegraphics[width=43mm,height=30mm]{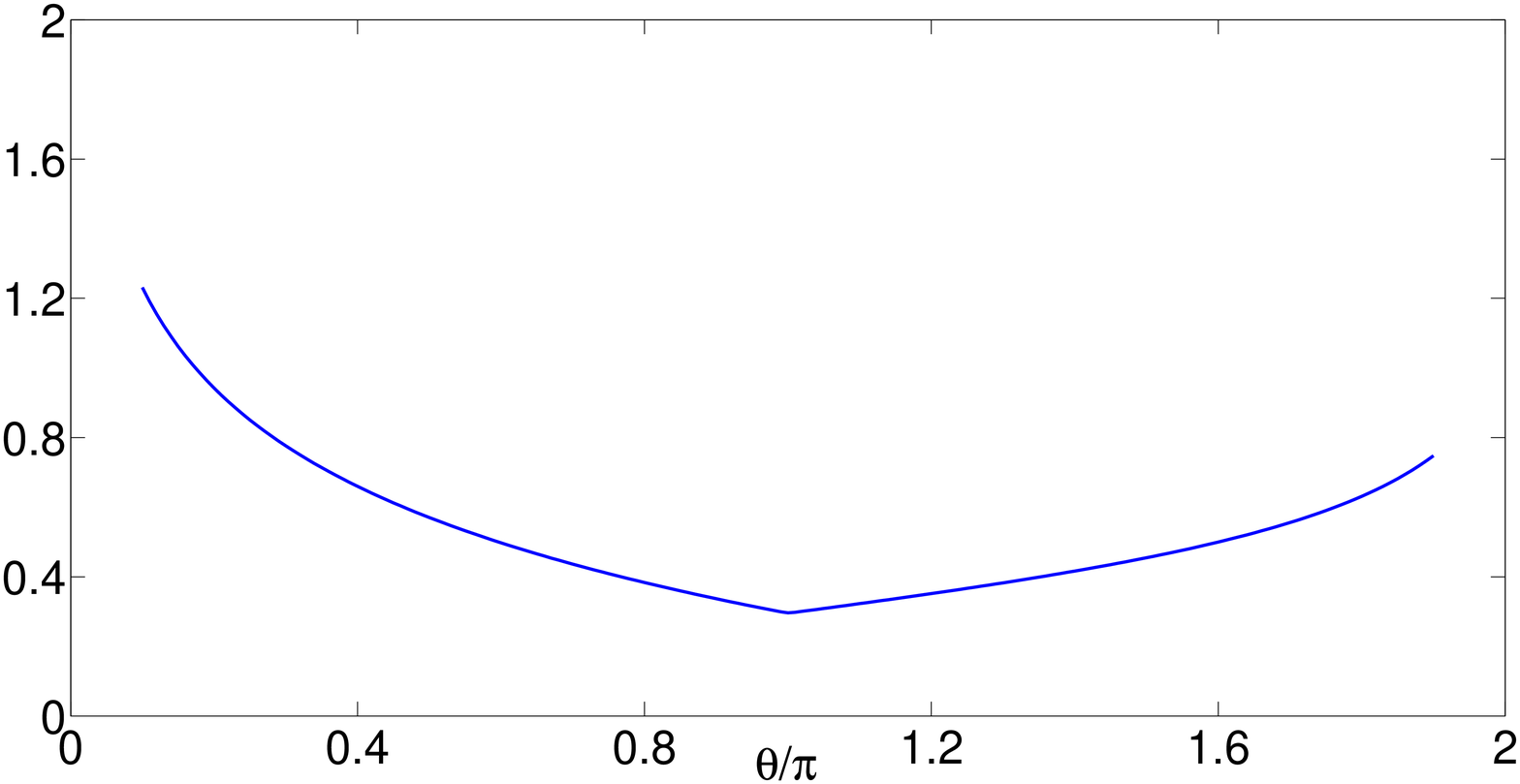}\hspace{-4mm}
\includegraphics[width=43mm,height=30mm]{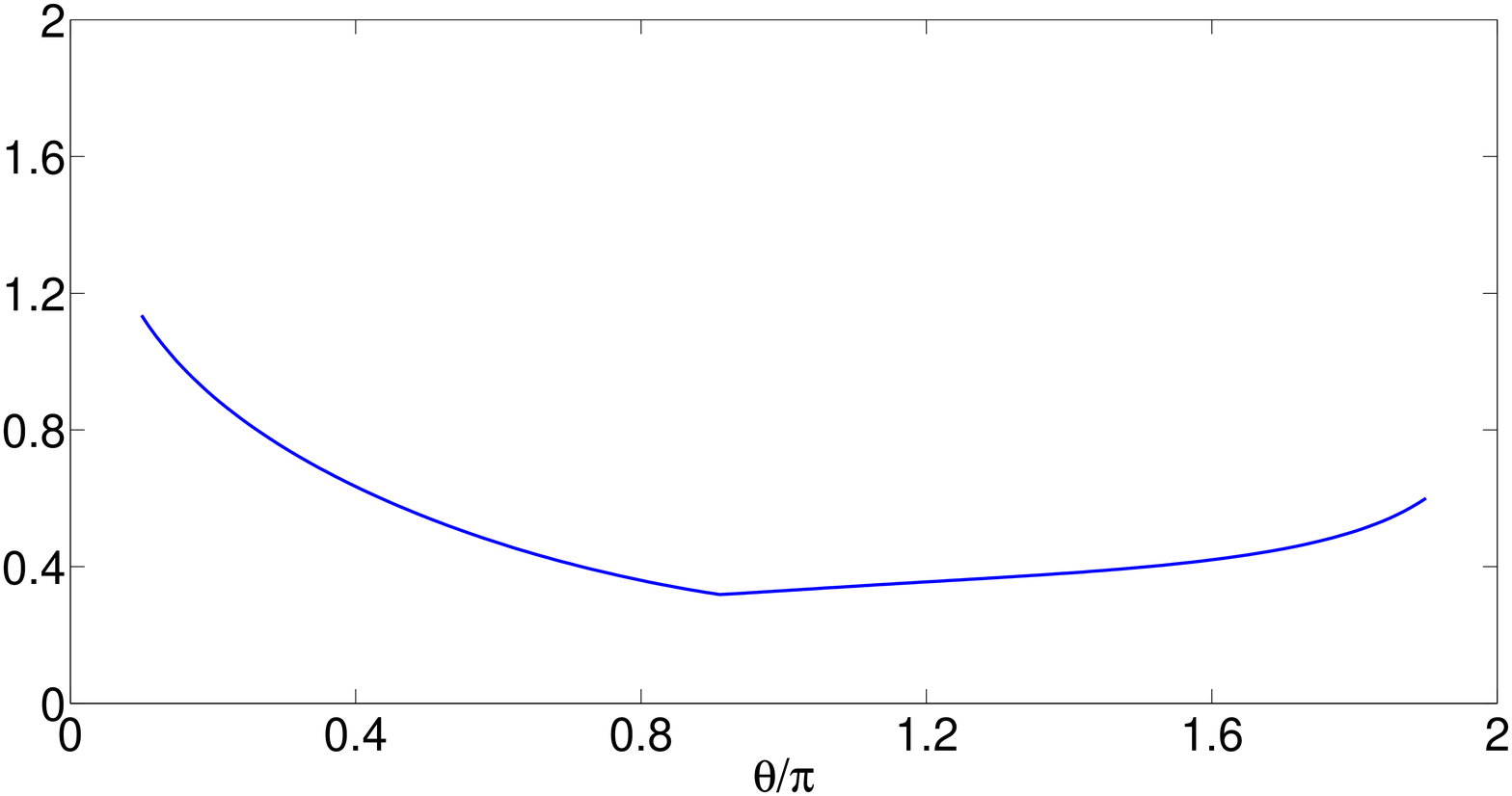}\\
 \includegraphics[width=43mm,height=30mm]{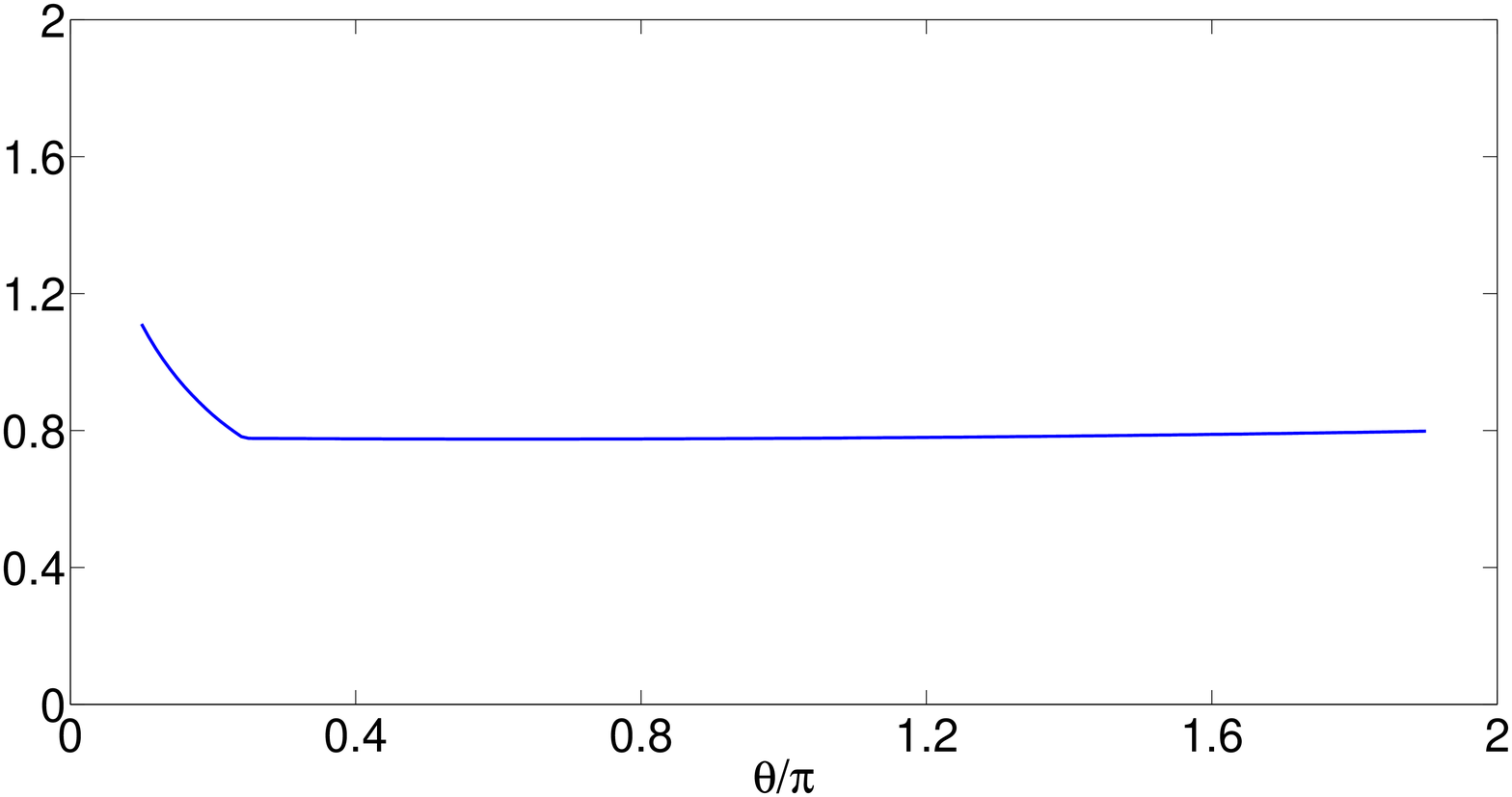}\hspace{-4mm}
 \includegraphics[width=43mm,height=30mm]{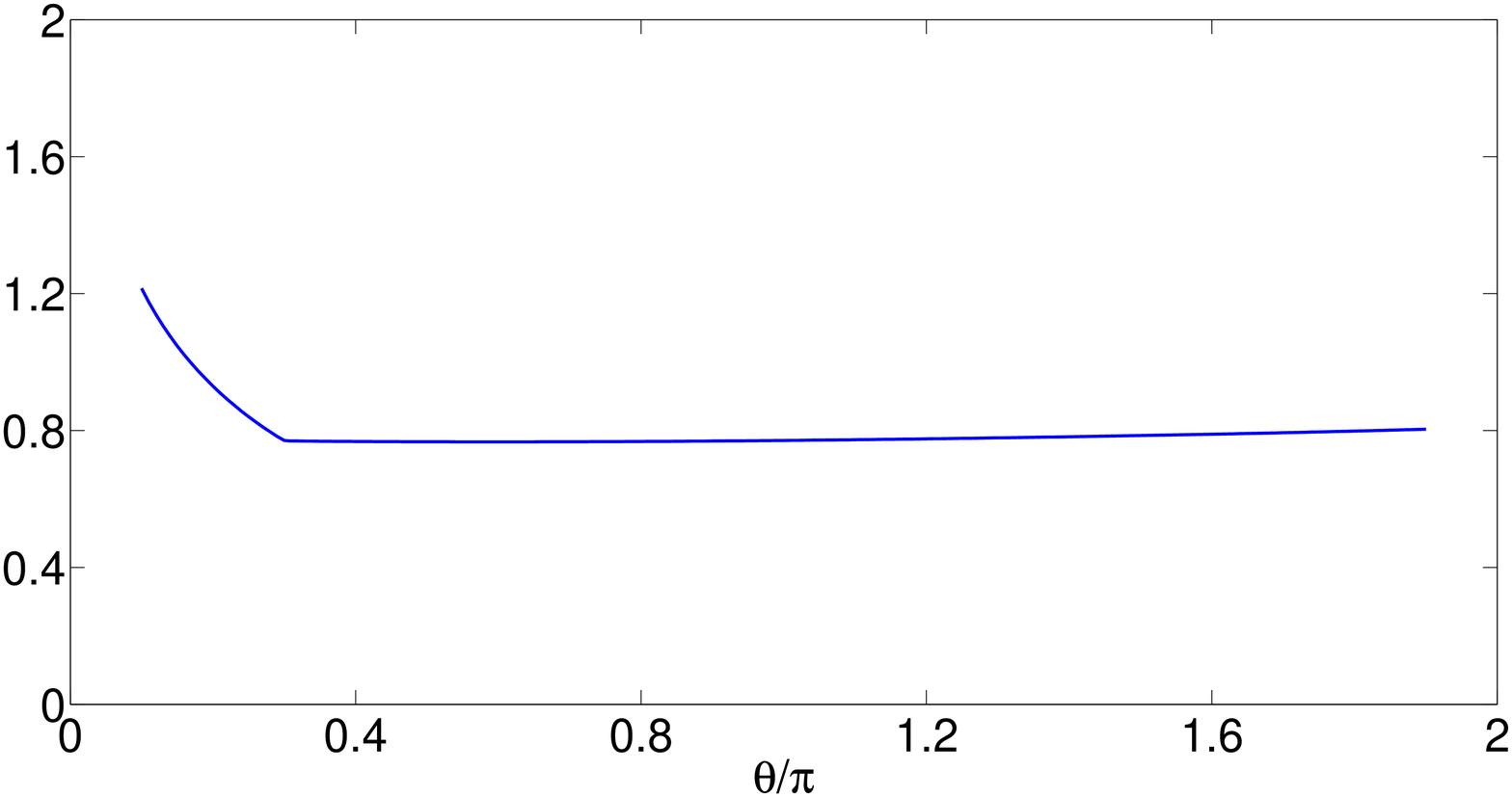}\hspace{-4mm}
 \includegraphics[width=43mm,height=30mm]{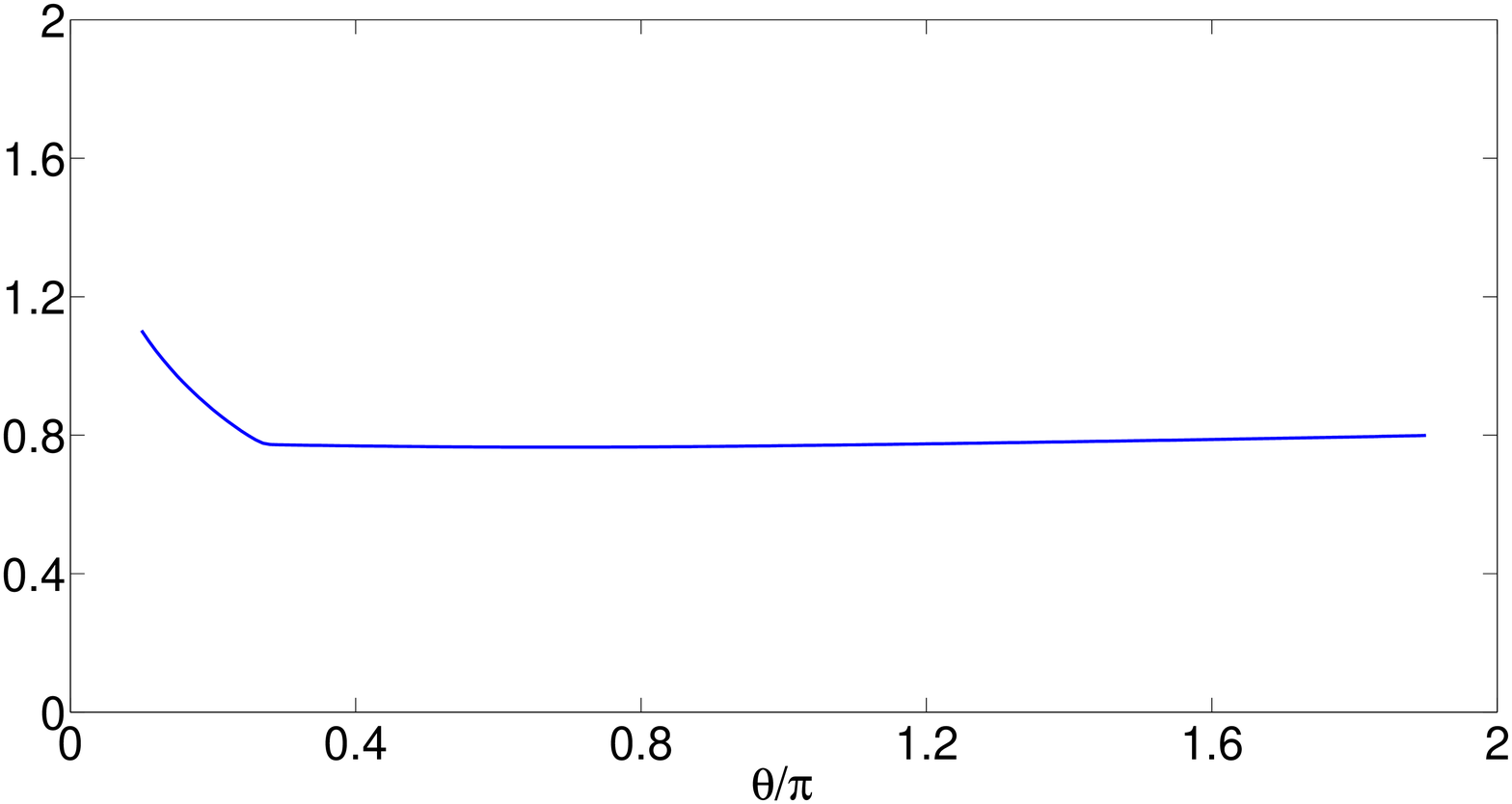}\\
 \includegraphics[width=43mm,height=30mm]{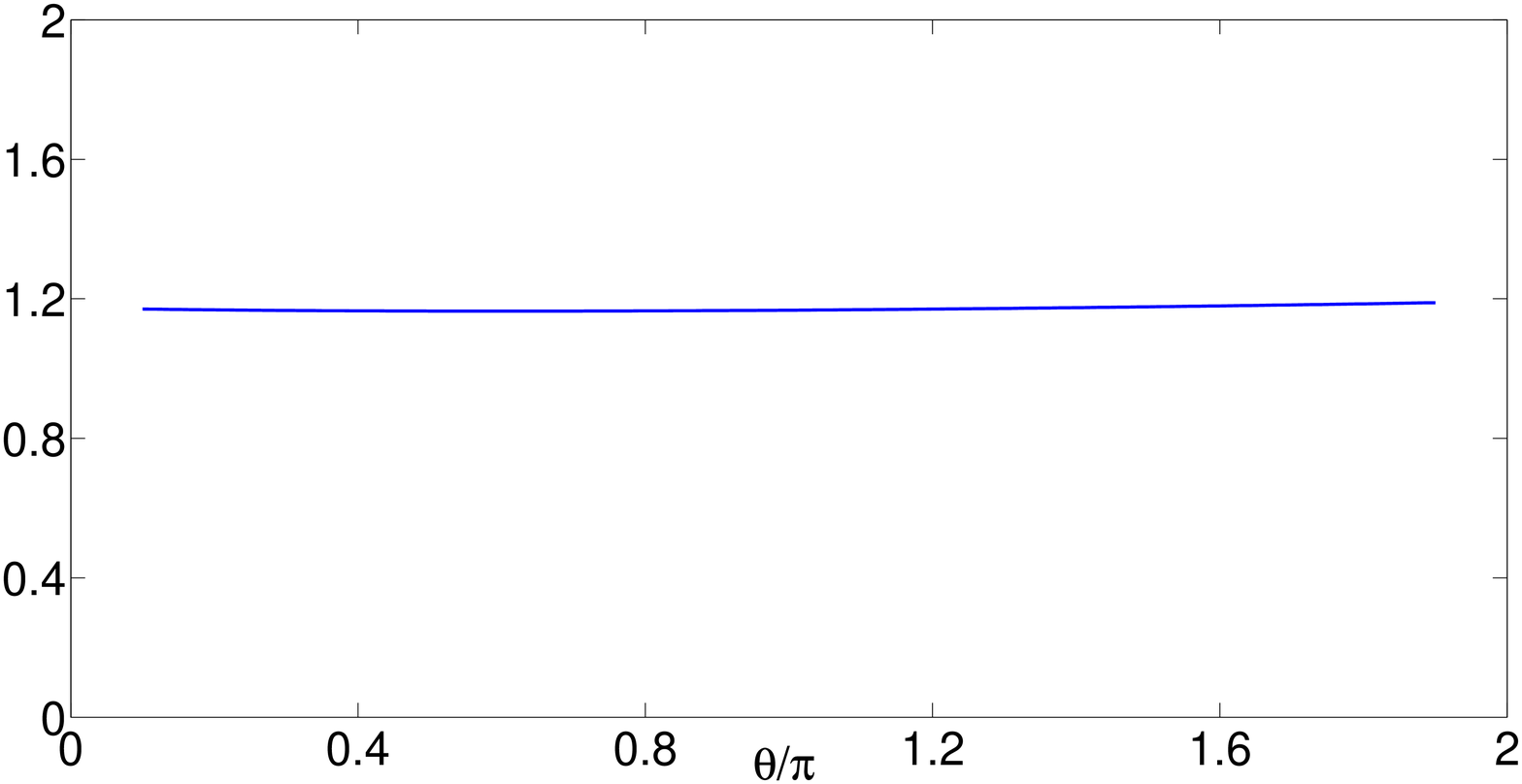}\hspace{-4mm}
 \includegraphics[width=43mm,height=30mm]{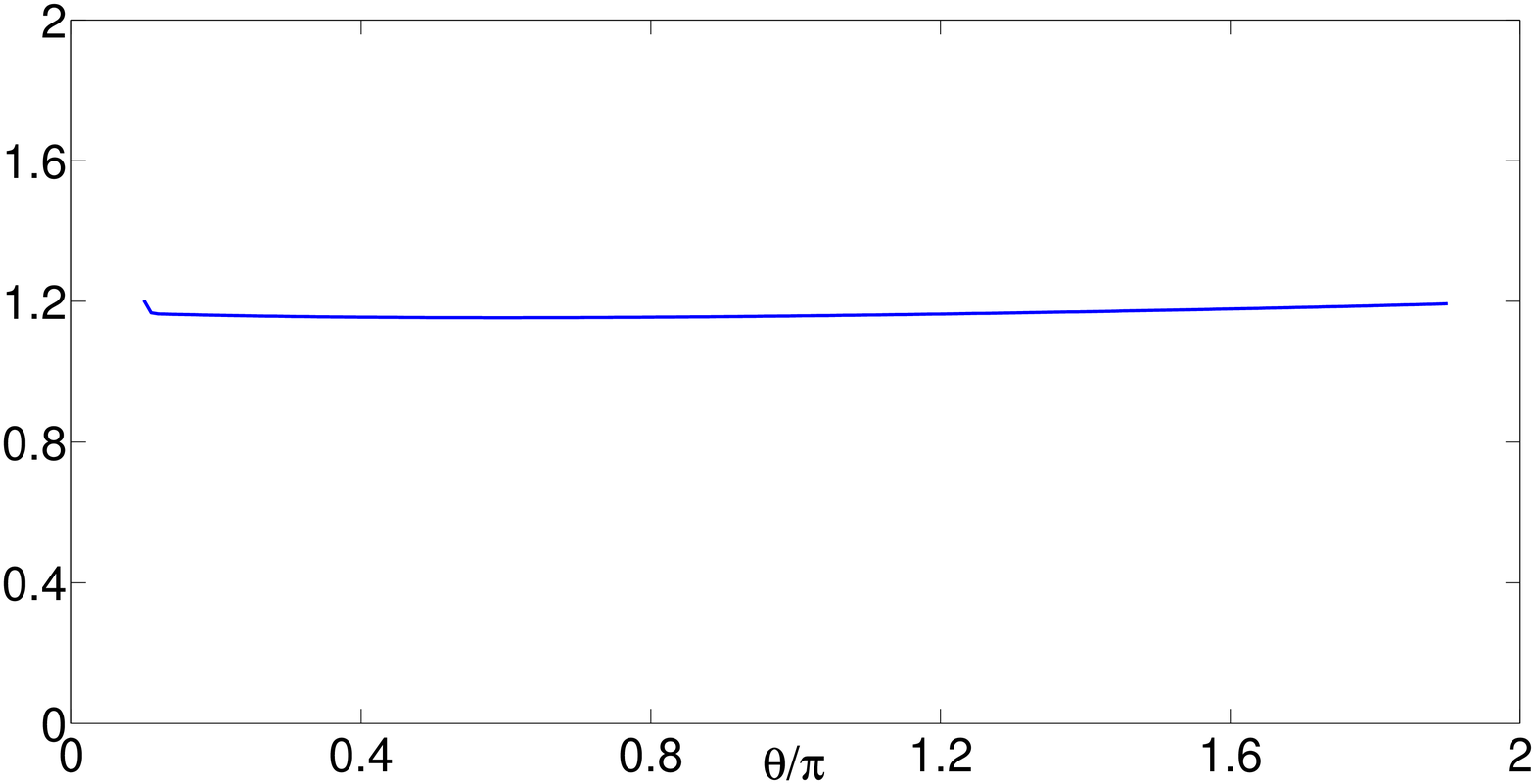}\hspace{-4mm}
 \includegraphics[width=43mm,height=30mm]{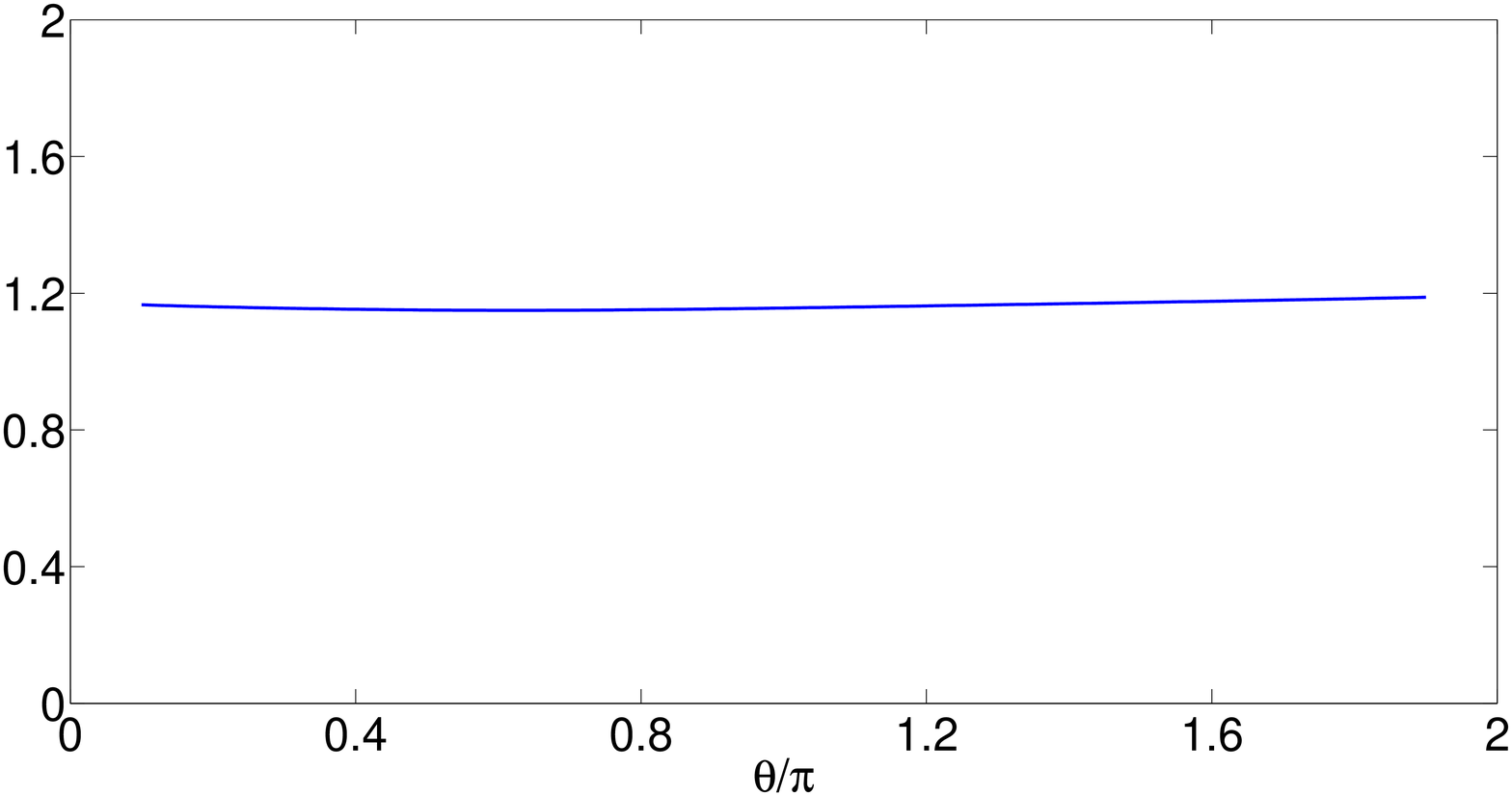}\\
  \caption{ Logarithm $\log_{10}$ of condition numbers  versus opening angles for Galerkin methods with $n=256$
  for various contours and spline spaces.
  Left: contour $\fL_1$; Middle: contour $\fL_2$; Right: contour
  $\fL_4$;
  First row: $d=0$; Second row: $d=1$; Third row: $d=2$}
 \label{cond_256}
\end{figure}

Now our numerical experiments can be described as follows. First, we
divide the interval $[0.1\pi,1.9\pi]$ by the points $\{\theta_k\}$,
where $\theta_k=\pi(0.1+k\times 0.01)$. For every family of the test
contours $\fL_j(\theta_k),\theta_k\in [0.1\pi,1.9\pi]$, $j=1,2,4$
consider the spline Galerkin methods with $n=256$ based on the
splines of order $0,1$ or $2$. Compute then the condition numbers of
the corresponding linear algebraic systems described by
\eqref{glk2}. Should there appear any point $\theta^*$ in the
vicinity of which the condition numbers become large, a neighborhood
of $\theta^*$ is refined by a smaller step $0.001\pi$, and condition
numbers are recalculated with $n$ changed to $512$. The outcome of
our computations is presented in Figure \ref{cond_256}. In all
cases, one can observe the absence of peaks in the graphs, which
means that the Galerkin methods under consideration do not have
"critical" angles in the interval $[0.1\pi,1.9\pi]$. In other words,
if the opening angles of all corners of the integration contour are
located in the interval  $[0.1\pi,1.9\pi]$, the spline Galerkin
methods based on the splines of degree $0,1$ or $2$ are always
stable. Another remarkable feature is that for the Galerkin methods
based on the splines of the same degree, all graphs are of the same
shape and for any $\theta \in [0.1\pi,1.9\pi]$ the corresponding
condition numbers are very close. This suggests a conjecture that
the condition numbers of the Galerkin methods possess certain
"locality" properties. They rather depend on the value of the
critical angles present than on the shape of the curves used.


 Note that all numerical experiments are performed in MATLAB
environment(version 7.9.0) and executed on an Acer Veriton M680
workstation equipped with a Intel Core i7 vPro 870 processor and 8GB
of RAM. These are time consuming computations and it took from one
to two weeks of computer work to obtain data for each graph in
Figure \ref{cond_256} .

\section{Conclusion}

In this work, necessary and sufficient conditions of the stability
of the spline Galerkin method for double layer potential equations
on simple piecewise smooth contours are established. The theoretical
results are verified by using curves with different number of corner
points and numerical results are in a good correlation with
theoretical studies. It turns out that the spline Galerkin methods
based on splines of degrees $0,1,2$ are always stable and the
convergence rates of the methods are comparable with other works.

\section{Acknowledgement} The authors would like to thank an anonymous
referee for constructive criticism and helpful suggestions.



 \end{document}